\documentclass[a4paper,12pt]{article}

\usepackage{geometry} 

\usepackage{titlesec}
\usepackage{float}
\titlespacing*{\section}
{0pt}{1.2ex plus 0.3ex minus 0.2ex}{0.8ex plus 0.2ex}

\titlespacing*{\subsection}
{0pt}{1.0ex plus 0.2ex minus 0.2ex}{0.6ex plus 0.2ex}

\usepackage{setspace} 

\usepackage{graphicx}
\usepackage{subcaption}
\usepackage{caption}

\usepackage{amsmath} 
\usepackage{amssymb}
\usepackage{mathtools}
\usepackage{amsthm}
\usepackage{enumitem}
\usepackage{url} 

\theoremstyle{remark}
\theoremstyle{plain}
\newtheorem{theorem}{Theorem}[section]   
\newtheorem{lemma}[theorem]{Lemma}

\usepackage{xcolor}
\usepackage{pifont}

\newcommand{\BubbleYes}{\textcolor{green!35!black}{\large\ding{51}}}
\newcommand{\BubbleNo}{\textcolor{red!85!black}{\large\ding{55}}}

\usepackage{graphicx}
\usepackage{subcaption}
\usepackage{placeins} 
\usepackage{pdflscape} 
\usepackage{amsthm}

\usepackage{booktabs} 

\usepackage[authoryear,round]{natbib} 
\usepackage[colorlinks=true, linkcolor=blue]{hyperref} 
    
\usepackage{comment}

\hypersetup{citecolor=blue,linkcolor=blue,urlcolor=blue} 

\graphicspath{ {./files/}} 

\newcommand{\tablepath}{files} 

\makeatletter
\newcommand{\includetable}[1]{%
  \@ifundefined{tablepath}{%
    \InputIfFileExists{#1}{}{}%
  }{%
    \InputIfFileExists{\tablepath/#1}{}{\InputIfFileExists{#1}{}{}}%
  }%
}
\makeatother

    \makeatother

\begin{document}
\begin{spacing}{1}
\clearpage
\thispagestyle{plain}

\thispagestyle{plain}

\onehalfspacing
    \begin{center}
        \LARGE
        \textbf{Double Local-to-Unity: Inference under Nearly Nonstationary Volatility }

        \Large
        Abir Sarkar\footnote{as4458@cornell.edu, 301 CIS Building,
Ithaca, New York 14853} and Martin T. Wells\footnote{mtw1@cornell.edu}  
        \vspace{0.25cm}
        
        Department of Statistics and Data Science, Cornell University 
        
        \vspace{0.25cm}
        
      \today
        
        \vspace{0.25cm}
    
    

        \begin{abstract}
            \footnotesize{\noindent 
            }

This article develops a moderate-deviation limit theory for autoregressive models with jointly persistent mean and volatility dynamics. The autoregressive root is allowed to drift toward unity more slowly than the classical \(1/n\) rate, while the stochastic-volatility persistence parameter also converges to one at a slower logarithmic rate. This double localization permits volatility to be nearly nonstationary, slowly evolving, and capable of generating divergent unconditional moments, features often observed in macro-financial data and asset-price bubble episodes. Under standard regularity conditions, we establish distributional limits for the ordinary least-squares estimator of the autoregressive coefficient under highly persistent, time-varying volatility. The effective normalization is governed by an average volatility scale, but the final least-squares limit laws are invariant to the nuisance parameters governing volatility specification. In the nearly stationary regime, where the autoregressive root approaches unity from below, the normalized OLS estimator is asymptotically normal. In the mildly explosive regime, where the root approaches unity from above, an OLS-based self-normalized statistic converges to a Cauchy limit. We also develop volatility-robust confidence intervals for autoregressive roots and local explosive-growth rates, providing feasible inference in nearly stationary and mildly explosive regimes. Simulations and empirical applications to historical bubble episodes, commodities, cryptocurrencies, and global real-estate markets show that our proposed procedure is less sensitive to volatility spikes and spurious explosive signals than standard homoskedastic bubble-detection approaches, while retaining power in sustained episodes of price exuberance.

        \end{abstract}
 \end{center}


\newpage

\section{Introduction}\label{section-intro}

Macroeconomic and financial time series often display two intertwined sources of dependence that challenge classical asymptotic theory: near-unit-root behavior in conditional means and persistent time-varying conditional variances. 
The former has been extensively analyzed in the unit-root and local-to-unity literature; see, for example, \cite{phillips1987towards} and \cite{chan1987asymptotic}. The second, equally important strand of the literature emphasizes that volatility is neither constant nor weakly persistent. Conditional heteroskedasticity can distort unit-root tests and alter the behavior of estimators and test statistics, motivating procedures that remain valid under broader classes of variance dynamics. \cite{kim1993unit} proved how conditional heteroskedasticity affects ADF-type tests in finite samples, highlighting the need for heteroskedasticity-robust formulations. Moderate deviations from a unit root towards explosiveness can be interpreted as asset price exuberance under homoskedasticity \cite{phillips2011explosive}, while a related literature on continuous-time asset pricing theory explored by \cite{JarrowProtterShimbo2010} emphasizes the pivotal role of volatility for the price process to carry an explosive bubble. This paper integrates these strands by developing inference on the autoregressive parameter under persistent stochastic volatility and examining its implications for bubble detection.

 \citet{phillips2007limit} studied the autoregressive root that approaches unity at moderate-deviation rate,
\[
y_t = \rho_n y_{t-1} + u_t, 
\qquad 
u_t = \sigma \varepsilon_t, 
\qquad 
\varepsilon_t \overset{\text{i.i.d.}}{\sim} (0,1),
\]
with $\rho_n = 1 + c/k_n$, and $k_n = o(n)$. This framework spans a continuum of cases, from the stationary, through local-to-unity dynamics, to mildly explosive regimes. In the nearly stationary moderate-deviation regime, Gaussian limit theory has been developed under several related parameterizations. \cite{giraitis2006uniform} obtain uniform normal limits for autoregressive coefficients satisfying $n(1-\rho_n)\to\infty$ under martingale-difference errors. \cite{park2003weak} considers the parameterization $\rho_n=1-m/n$, with $m,n\to\infty$, and derives asymptotic normality of the sample correlation coefficient at rate $n/\sqrt{m}$ when $m=o(n)$. In the same moderate-deviation framework, \cite{phillips2007limit} show that when the root approaches unity from below at rate $1/k_n$, the least-squares estimator is asymptotically normal under the normalization $\sqrt{n k_n}$.

In the mildly explosive regime, $\rho_n=1+c/k_n$ with $c>0$, \cite{phillips2007limit} proved that $\widehat{\rho}_n-\rho_n$ converges at rate $k_n\rho_n^{\,n}$ and that the associated serial correlation statistic has a Cauchy limit. This result connects the moderate-deviation framework to the classical theory of explosive autoregressions: for fixed $|\rho|>1$, \cite{white1958limiting} obtains a Cauchy limit for the Gaussian AR(1), while \cite{basawa1984asymptotic} develops related conditional limit theory for nonergodic stochastic processes. The intuition is that explosive growth in the regressor dominates weak dependence in the innovations. Accordingly, \cite{phillips2007limitb} extend the Cauchy limit to linear-process errors
\[
u_t=\sum_{j=0}^{\infty}c_j e_{t-j},
\qquad
\mathbb{E}[e_t]=0,\quad
\mathbb{E}[e_t^2]<\infty,\quad
\sum_{j=0}^{\infty}j|c_j|<\infty.
\]
Relatedly, \cite{mikusheva2012one} studies testing and confidence set construction for functions of autoregressive coefficients in persistent time-series models.

Time-varying volatility is central to estimation and inference in persistent autoregressive models. \cite{cavaliere2009heteroskedastic} develop asymptotic theory for Dickey--Fuller statistics under nonstationary stochastic volatility and propose wild-bootstrap inference, while the near-unit persistence often observed in SV and GARCH models has motivated nearly integrated volatility models such as IGARCH and GARCH near the unit-root boundary \citep{nelson1990stationarity}. The asymptotic tools for these settings build on martingale and linear-process methods; in particular, \cite{phillips1992asymptotics} provide a unified treatment of strong laws, central limit theory, and invariance principles for dependent heteroskedastic innovations. More generally, \citet{dahlhaus1997fitting} studies nonstationary time series with time-varying second moments via evolutionary spectra, while \citet{roy2026asymptotic} provides minimum-divergence tools for inference under misspecification. Recent work also incorporates volatility directly into inference on persistent and explosive dynamics, including de-volatilized recursive procedures based on high-frequency volatility measures \citep{boswijk2024testing}. Together, this literature points to a central econometric issue: inference on persistent autoregressive coefficients can depend critically on the volatility environment, especially when volatility is nearly nonstationary.

\subsection{Econometric Implications}

Moderate deviations from a unit root indicate mildly explosive dynamics in the underlying process.  In financial markets, such departures are commonly interpreted as price exuberance, i.e., an asset price bubble. The asymptotic foundations for mildly explosive autoregressions were developed by \citet{phillips2009limit}, while \citet{phillips2011explosive} introduced the PWY recursive testing procedure for detecting explosiveness and dating the origination and collapse of exuberance, with an application to the 1990s Nasdaq episode. This framework was subsequently extended by \citet{phillips2015testing} to real-time dating (PSY) algorithms that allow for multiple bubble episodes. Related methodological developments and applications include \citet{harvey2016tests}, \citet{astill2018real}, \citet{whitehouse2019explosive}, \citet{evgenidis2020lean}, and \citet{whitehouse2023real}.

The empirical reach of this literature is now broad. Recursive explosive-root methods have been used to study housing-market exuberance in Australia, China, Hong Kong, and the United States \citep{shi2016dating,deng2017did,yiu2013detecting,shi2017speculative}, cryptocurrency bubbles \citep{corbet2018datestamping}, common bubble behavior in high-dimensional financial systems \citep{chen2023common}, and the transmission of cryptocurrency bubbles to systemic risk in global energy companies \citep{ji2022cryptocurrency}. More recently, \citet{basele2025speculative} apply these ideas to speculative episodes among the ``Magnificent Seven'' technology firms during the AI boom. 

A key limitation of much of the bubble-detection literature is its baseline homoskedasticity assumption, which is difficult to reconcile with asset-price data. Volatility clustering, regime shifts, and leverage effects are well documented in financial markets \citep{engle1982autoregressive,engle1986modelling}, and realized volatility often rises sharply during periods of exuberance and stress \citep{greenwood2019bubbles,shephard2010realising,patton2015good}. Recent work has therefore begun to incorporate volatility more explicitly into bubble analysis, including continuous-time no-arbitrage restrictions linking explosive prices to volatility conditions \citep{JarrowProtterShimbo2010}, option-implied methods for measuring bubble magnitudes under stochastic volatility \citep{jarrow2021inferring,fusari2025testing}, quadratic-variation approaches to explosive behavior \citep{jarrow2024study}, cryptocurrency bubble-dating procedures benchmarked against GARCH specifications \citep{choi2020testing,montanino2025bubble}, and  bubble-dating applications to recent technology and semiconductor firms \citep{sarkar2026there}. These developments motivate an autoregressive theory of economic bubble detection that remains valid when volatility is persistent, time varying, and potentially close to nonstationary.


\subsection{Our Contributions}\label{sec:our_conts}

We study the AR(1) model with a near--nonstationary stochastic--volatility recursion
\begin{equation}\label{eq:first model}
\begin{split}
y_t &= \rho_n y_{t-1} + u_t, \quad u_t = \sigma_t \varepsilon_t, \quad 
      \varepsilon_t \overset{\text{i.i.d.}}{\sim} (0,1), \\
\log \sigma_t^2 &= \phi_n \log \sigma_{t-1}^2 + \eta_t, \quad
\eta_t \overset{\text{i.i.d.}}{\sim}  \mathcal{N}(0,\alpha^{2}), \quad
\eta_t \perp \varepsilon_t, \quad \alpha > 0 .
\end{split}
\end{equation}
The AR(1) log-volatility specification is canonical in time-series econometrics and asset pricing, offering a parsimonious representation of volatility clustering and persistence \citep{taylor2008modelling,harvey1994multivariate}. We impose no restrictive distributional assumption on the innovation $\varepsilon_t$, in line with the moderate-deviation autoregression literature; it is enough that $\mathbb{E}[\varepsilon_t]=0$, $\mathbb{E}[\varepsilon_t^2]=1$, and $\mathbb{E}[\varepsilon_t^4]<\infty$. The Gaussian assumption on $\eta_t$ is used only to obtain a tractable lognormal stochastic-volatility scale and is standard in AR(1) stochastic-volatility models. Both persistence parameters are allowed to drift with the sample size. The mean coefficient satisfies $\rho_n=1\pm c/k_n$ with $k_n\to\infty$ and $k_n=o(n)$, while the volatility coefficient is allowed to follow an ultra-persistent path with $\phi_n\uparrow 1$. This generates a double-local-unit environment in which both the conditional mean and the conditional variance approach nonstationarity. Such persistence is consistent with empirical evidence from highly persistent stochastic-volatility and GARCH-type processes with persistence factor close to, but below, one \citep{bollerslev1994arch,baillie1996fractionally}.

The main novelty of our paper is that the familiar moderate-deviation limit theory for autoregressive estimation is preserved under highly persistent volatility. Although the volatility process is nearly nonstationary and its unconditional moments may diverge, the final least-squares limits retain the same Gaussian form in the nearly stationary regime and the same Cauchy form in the mildly explosive regime. The difficulty is that this robustness is not automatic,  volatility changes the intermediate normalization, creates non-negligible cross-moment dependence, and makes the martingale variance random.

We identify the stochastic-volatility normalizations that make this
robustness possible. In the nearly stationary regime, the relevant scale is the
endogenous average volatility
\begin{equation}\label{eq:mn_def_intro}
m_n 
:= \frac{1}{n}\sum_{t=1}^{n}\mathbb{E}[\sigma_t^2].
\end{equation}


The main technical challenge is that persistent volatility changes the intermediate normalization, creates non-negligible dependence in volatility cross moments, and makes the predictable variation of the OLS score random. We establish the lognormal moment bounds, quadratic-variation law of large numbers, and volatility-weighted martingale central limit theorem needed to control these effects. The resulting robustness mechanism is important. Once the correct scale \(m_n\) is used, the volatility contribution enters the sample covariance, sample variance, and martingale variance in the same way, and therefore cancels from the final Gaussian limit, yielding nuisance-parameter free distribution even though volatility is highly persistent and time varying.

We develop the corresponding mildly explosive theory under persistent
stochastic volatility. In this regime the relevant long-run scale is
\begin{equation}\label{eq:elln_def_intro}
    \ell_n:=\exp\left\{\frac{\alpha^2}{2(1-\phi_n^2)}\right\}.
\end{equation}

After normalization by \(\ell_n\), the covariance and quadratic-variation terms
recover the product-normal and quadratic-normal structure of the homoskedastic
mildly explosive theory. The key difficulty is to show that volatility-induced remainder terms and triangular double sums involving \(u_j u_t\) are negligible relative to the explosive rate \(\rho_n^n\). Once this is established, the self-normalized OLS statistic retains the standard Cauchy limit,  as in classical results such as \cite{anderson1959asymptotic} and \cite{white1958limiting},
despite the presence of nearly nonstationary stochastic volatility.

Building on the seminal analysis of \cite{phillips2007limit}, these results extend moderate-deviation asymptotics from homoskedastic autoregressions to a setting in which both the conditional mean and the conditional variance drift toward nonstationarity. The contribution is not merely to allow heteroskedastic errors, but to handle a volatility mechanism whose moments and temporal dependence evolve with the sample size. This yields feasible inference for persistent and mildly explosive autoregressive dynamics without estimating the nuisance parameters governing volatility specifications, and broadens the applicability of moderate-deviation methods to empirical settings with unstable volatility. In our real-data applications, this volatility-robust perspective produces economically more plausible and statistically more stable bubble classifications than standard homoskedastic bubble-detection approaches.


The remainder of the paper is organized as follows.  In Section~\ref{sec:notations}, we define all the mathematical notations used throughout the paper. Section~\ref{section-methods} discusses the key assumptions and a brief proof idea involving the double local-to-unity asymptotics. Section~\ref{sec:section3} and~\ref{sec:section4} develop the properties involving distribution of the autoregressive coefficient and explosive growth parameters under different persistent mean and volatility regimes. Section~\ref{sec:simulation} provides empirical evidence through simulations and real-data studies that justify the theoretical results with statistical analyses, and we conclude the paper in Section~\ref{sec:conclusion} with a summary of novel contributions and a discussion of potential directions for future research. All of the proofs of the main results can be found in the \hyperref[sec:appendix]{Appendix}, and the code for empirical analysis on GitHub.

\subsection{Notations}\label{sec:notations}
All random quantities are defined on a common probability space \((\Omega,\mathcal{F},\mathbb{P})\). We write \(\mathbb{P}\) for probability and \(\mathbb{E}\) for expectation, with conditional expectations taken relative to the indicated \(\sigma\)-algebras. The filtration
\(
\mathcal{F}_{t}
:=
\sigma\big(y_0,\varepsilon_1,\ldots,\varepsilon_t,\eta_1,\ldots,\eta_t\big)
\)
denotes the information available up to time \(t\). For deterministic sequences, \(a_n=o(b_n)\) means \(a_n/b_n\to0\), \(a_n=\mathcal{O}(b_n)\) means \(|a_n/b_n|\) is bounded, \(a_n\asymp b_n\) means the ratio is bounded away from zero and infinity, \(a_n\sim b_n\) means \(a_n/b_n\to1\), and \(a_n\gg b_n\) means \(a_n/b_n\to\infty\). For random sequences, \(X_n=o_p(1)\) denotes convergence to zero in probability, \(X_n=\mathcal{O}_p(1)\) denotes boundedness in probability, \(X_n\xrightarrow{p}X\) denotes convergence in probability, and \(X_n\xRightarrow{d}X\) denotes convergence in distribution. We also write \(X_n\lesssim Y_n\) when \(X_n\le C Y_n\) for some constant \(C>0\) independent of \(n\). Finally, \(\mathbf{1}\{\cdot\}\) is the indicator function, \(\lfloor\cdot\rfloor\) is the floor function,  \(\mathrm{poly}(a_n)\) denotes any polynomial in \(a_n\), \(\mathcal{N}(\mu,\sigma^2)\) denotes a normal law, and \(\mathcal{C}(0,1)\) denotes the standard Cauchy law with density \(f(x)=\{\pi(1+x^2)\}^{-1}\).

\section{Deviations from unity and key assumptions}\label{section-methods}

We study model~\eqref{eq:first model} when both the autoregressive root and the volatility persistence parameter drift toward unity. The mean coefficient satisfies the moderate-deviation parameterization
\begin{equation}\label{eq:mean_eq}
\rho_n = 1 \pm \frac{c}{k_n},
\qquad 
c>0,\qquad 
k_n\to\infty,\qquad 
k_n=o(n),
\end{equation}
so that the root approaches unity more slowly than in the classical $1/n$ local-to-unity framework. The volatility coefficient is specified as
\begin{equation}\label{eq:vol_eq}
\phi_n = 1-\frac{d}{\log r_n},
\qquad 
d>0,
\end{equation}
where \(r_n>1\) and, in the nearly nonstationary volatility case, \(r_n\to\infty\). One must note, constant $\phi_n$ with $|\phi_n| < 1$ is allowed in this setting. Moreover,~\eqref{eq:vol_eq} allows the log-volatility process highly persistent while remaining mean reverting. Thus, the model can potentially have two drifting components. We now impose the following assumptions.

\newenvironment{assumptions}{%
\par\noindent\textbf{Assumptions.}%
\begin{enumerate}%
\renewcommand\labelenumi{A\arabic{enumi}.}%
\renewcommand\theenumi{A\arabic{enumi}.}%
}{%
\end{enumerate}%
}

\begin{assumptions}
\item \label{assm:A1} Nearly stationary mean: 
\[
\rho_n = 1-\frac{c}{k_n},
\qquad 
\phi_n = 1-\frac{d}{\log r_n},
\qquad 
c,d>0.
\]

\item \label{assm:A2} Mildly explosive mean:
\[
\rho_n = 1+\frac{c}{k_n},
\qquad 
\phi_n = 1-\frac{d}{\log r_n},
\qquad 
c,d>0.
\]

\item \label{assm:A3} Rate separation:
\[
r_n^q=o(k_n),
\qquad
k_n r_n^q=o(n),
\qquad \text{for every fixed } q>0.
\]
\end{assumptions}

Assumption~\hyperref[assm:A1]{A1} is used for the nearly stationary theory, while Assumption~\hyperref[assm:A2]{A2} is used for the mildly explosive theory. Assumption~\hyperref[assm:A3]{A3} imposes the rate separation between the mean-root localization and the volatility persistence. It requires every fixed power of \(r_n\) to be asymptotically dominated by both \(k_n\) and \(n/k_n\). Thus, the volatility process is allowed to be highly persistent without dominating the mean dynamics. Intuitively, this places \(\rho_n\) in a moderate-deviation neighborhood of unity---persistent enough to require nonstandard normalization, but far enough from the exact unit-root boundary to yield nondegenerate limits.

To analyze the model~\eqref{eq:first model} under these assumptions, we focus on the least-squares estimator of the autoregressive coefficient. Although \(\rho_n\) drifts toward unity, it does so at a moderately slow rate, so the process lies in a moderate-deviation region between the classical fixed root autoregression and local-to-unity cases, with a separate mildly explosive limit when the root approaches unity from above. The object of interest is
\begin{equation}\label{eq:rho_eq}
\widehat{\rho}_n
=
\frac{\sum_{t=1}^{n}y_{t-1}y_t}
     {\sum_{t=1}^{n}y_{t-1}^{2}},
\qquad
\widehat{\rho}_n-\rho_n
=
\frac{\sum_{t=1}^{n}y_{t-1}u_t}
     {\sum_{t=1}^{n}y_{t-1}^{2}}.
\end{equation}

We now briefly outline the proof strategy, emphasizing how the preceding robustness mechanisms enter the asymptotic arguments. The proof proceeds by studying the two components in~\eqref{eq:rho_eq}: the quadratic denominator \(\sum_{t=1}^{n}y_{t-1}^{2}\) and the martingale score \(\sum_{t=1}^{n}y_{t-1}u_t\). After controlling the sample covariance term, quadratic variation in the denominator satisfies
\begin{equation}\label{eq:denom_idea_stationary}
\frac{1}{n k_n m_n}
\sum_{t=1}^{n}y_{t-1}^2
\xrightarrow{p}
\frac{1}{2c}.
\end{equation}
Thus, the denominator retains the familiar moderate-deviation form, but only after the volatility scale induced by the nearly nonstationary variance process is incorporated. The numerator is handled through a martingale CLT. Although \(y_{t-1}u_t\) remains a martingale difference, its conditional variance is volatility dependent:
\begin{equation}\label{eq:cond_var_idea}
\mathbb{E}\!\left[u_t^2\mid \mathcal{F}_{t-1}\right]
=
e^{\alpha^2/2}\sigma_{t-1}^{2\phi_n}.
\end{equation}
Consequently, the predictable variation involves a volatility-weighted quadratic form, and the key auxiliary convergence is
\begin{equation}\label{eq:weighted_lln_idea}
\frac{1}{n k_n m_n^2}
\sum_{t=1}^{n}y_t^2\sigma_{t-1}^{2\phi_n}
\xrightarrow{p}
\frac{e^{-\alpha^2/2}}{2c}.
\end{equation}
Together with a Lindeberg condition, this gives the Gaussian limit for the centered least-squares estimator. The important feature is that the same volatility scale enters the numerator and denominator; hence it cancels from the final limit.

Under Assumption~\hyperref[assm:A2]{A2}, the same decomposition is used, but the leading terms are governed by the terminal behavior of explosive autoregression. After normalization by \(\ell_n\), the covariance and the quadratic variation terms can be written, up to negligible remainders, in terms of weighted innovation sums that converge jointly to independent Gaussian limits. The covariance term in numerator has a product-normal limit, while the quadratic variation term in denominator has a quadratic-normal limit, and their ratio yields the Cauchy limit  as in the homoskedastic mildly explosive theory.

A central implication is that inference on the autoregressive root remains robust to the nuisance volatility scale generated by nearly nonstationary stochastic volatility. Since the Gaussian and Cauchy limits coincide with their homoskedastic counterparts, the confidence-interval methods for \(\rho_n\) and local explosive-growth rate $\gamma_n$ developed in \citet{phillips2023estimation} remain applicable in this setting. Thus, empirical tools for bubble-growth estimation remain feasible under time-varying volatility and deliver diagnostics that are less sensitive to volatility spikes. In practice, this yields more reliable uncertainty quantification and a more robust detection of exuberance near the explosive boundary.

\section{Double local-to-unity in nearly stationary regime}\label{sec:section3}

This section analyzes the large-sample behavior of the autoregressive coefficient and the growth parameter under model \eqref{eq:first model} with Assumptions~\hyperref[assm:A1]{A1} and~\hyperref[assm:A3]{A3}. We will stick to this set of assumptions throughout this section.

We initialize the process at \(y_0\) with the condition $\mathbb{E}[y_0^4]=o(k_n^2)$ and $\sigma_0 =1$, both independent of the $\sigma$-field generated by $\{(\varepsilon_t, \eta_t): t \geq 1\}$. The fourth-moment condition keeps the contribution of the initial state small enough relative to the moderate-deviation scale, so that the limiting behavior is driven by the accumulated innovations rather than by just $y_0$. The sequences $\{\varepsilon_t\}$ and $\{\eta_t\}$ are independent of each other. Because both $\rho_n$ and $\phi_n$ depend on $n$, the model defines triangular arrays $\{y_{nt}: 1 \leq t \leq n\}$. For notational simplicity we write $y_t$, suppressing the $n$ index. We first show that, under an appropriate normalization tailored to the local-to-unity rate $k_n$, asymptotic behavior of $\widehat{\rho}_n$ in~\eqref{eq:rho_eq} parallels that of the strictly stationary case $\rho < 1$, providing the key scale for inference on $\rho_n$ even in the presence of highly persistent volatility. We need to define a few mathematical notations.

\noindent Let \(z_t:=\log\sigma_t^2\). Iterating the AR(1) volatility recursion gives
\begin{equation}\label{eq:expression_z}
z_t
=
\phi_n^t z_0+\sum_{j=0}^{t-1}\phi_n^j\eta_{t-j},
\qquad
\sigma_t^2
=
\exp\!\left(\phi_n^t z_0+\sum_{j=0}^{t-1}\phi_n^j\eta_{t-j}\right).
\end{equation}
Define
\begin{equation}\label{eq:x_m_A_def}
A_t
:=
\frac{1-\phi_n^{2t}}{2(1-\phi_n^2)},
\qquad
x_t
:=
\mathbb{E}[\sigma_t^2]
=
\exp(\alpha^2 A_t),
\qquad
m_n
=
\frac{1}{n}\sum_{t=1}^{n}x_t .
\end{equation}

\begin{lemma}\label{lemma:main_lemma_1}
    For the model in \eqref{eq:first model} with Assumptions~\hyperref[assm:A1]{A1} and \hyperref[assm:A3]{A3}, uniformly over \(1\le t\le n\),
    \begin{enumerate}
\renewcommand{\labelenumi}{(\alph{enumi})}
   \item $ \mathbb{E}[\sigma_t^2] \leq e^{\alpha^2A_n} \;\lesssim\; r_n^{\alpha^2/(4d)}$,
   \item $ \label{lemma:b} \mathbb{E}[\sigma_t^4] \leq e^{4\alpha^2 A_n}  \;\lesssim\; r_n^{\alpha^2/d}$,
   \item  \label{lemma:c} $\frac{1}{n^2} \sum \sum_{t \ne s} \mathbb{E}[\sigma_t^2 \sigma_s^2] \lesssim e^{2\alpha^2 A_n}$.
    \end{enumerate}
\end{lemma}
The first two bounds control the second and fourth-order moments of the volatility process, preventing growth beyond the required rate and ensuring tightness of variance-driven quantities. The cross-moment bound controls cross moments over time, limiting temporal dependence in double sums and stabilizing the variance of aggregates. This bound is crucial in particular, as it  controls temporal dependence in volatility-driven double sums, which is the key reason for the asymptotic distribution to stabilize. Together, these bounds provide the concentration, rate control, and negligibility of cross terms needed for the proof.

\begin{lemma}\label{lemma:mean_lemma}
    Under model \eqref{eq:first model} with Assumptions~\hyperref[assm:A1]{A1} and \hyperref[assm:A3]{A3},
    $m_n = e^{\alpha^2 A_n} \,(1+ o(1))$.
    \end{lemma}

This result is central in establishing the consistency of \(\widehat{\rho}_n\). The quantity \(m_n\) is the average unconditional innovation variance and therefore represents the effective volatility scale of the sample when \(\{\sigma_t^2\}\) is time varying and persistent. In this sense, \(m_n\) provides the natural normalization under stochastic volatility: after rescaling by \(m_n\), the quadratic terms have the same nondegenerate order as in the homoskedastic moderate-deviation theory. In the usual double-logarithmic case, \(r_n=\log n\), one can show \(m_n\) grows only polylogarithmically in \(n\), so it modifies the scale without changing the fundamental \(\sqrt{n k_n}\) rate obtained in \cite{phillips2007limit}.  Economically, \(m_n\) summarizes average volatility risk over the sample. A larger \(m_n\) reflects prolonged high-uncertainty episodes, implying wider confidence intervals and more data needed to detect mean reversion. This connects to continuous-time asset-pricing work linking bubble behavior to restrictions on volatility and mean reversion \citep{MijatovicUrusov2012,JarrowProtterShimbo2010}.

Now, we should note that, the condition $\mathbb{E}[y_0^4]=o(k_n^2)$ implies \(y_0=o_p(\sqrt{k_n})\) and makes the initial condition asymptotically negligible under the normalizations used below.

\begin{lemma}\label{lem:main_lemma_2} Let $y_0=o_p\left(k_n^{1/2}\right)$. For the model in \eqref{eq:first model}, under the Assumptions~\hyperref[assm:A1]{A1} and~\hyperref[assm:A3]{A3},
\begin{enumerate}
\renewcommand{\labelenumi}{(\alph{enumi})}
\item\label{lem:lemma_main_1a}
If $k_n \,r_n^{\alpha^2/(4d)}=o(n)$, then $\frac{y_n^2}{n}\xrightarrow{p}0$.
\item \label{lem:lemma_main_1b}If $ k_n\,r_n^{\frac{\alpha^2}{d}} = o(n)$, then $\displaystyle \frac{1}{n} \sum_{t=1}^{n} y_{t-1} u_t \xrightarrow{p} 0.$
\end{enumerate}
\end{lemma}

\noindent We will leverage this lemma later in the proof. First, we aim to show that, \begin{align}\label{eq:equation6}
\frac{1}{n k_n  m_n} \sum_{t=1}^{n} y_{t-1}^2 
\xrightarrow{p} \frac{1}{2c}. 
\end{align}
From the identity
\(
y_t^2 - \rho_n^2 y_{t-1}^2 
= u_t^2 + 2\rho_n y_{t-1}u_t,
\)
we obtain, upon summing over $t=1,\dots,n$ and dividing by $n$,
\begin{align}\label{eq:equ5_main}
(1 - \rho_n^2)\cdot \frac{1}{n}\sum_{t=1}^{n} y_{t-1}^2
&= \frac{1}{n}(y_0^2 - y_n^2)
+ \frac{1}{n}\sum_{t=1}^{n} u_t^2
+ 2\rho_n\frac{1}{n}\sum_{t=1}^{n} y_{t-1} u_t \notag \\[6pt]
&= \frac{1}{n} \sum_{t=1}^{n} u_t^2 + o_p(1),
\end{align}
where we used $y_0 = o_p(k_n^{1/2})$, $y_n^2/n \xrightarrow{p} 0$ (by Lemma \hyperref[lem:lemma_main_1a]{3.3 (a)}), and 
$\frac{1}{n}\sum_{t=1}^{n} y_{t-1}u_t = o_p(1)$ (by Lemma \hyperref[lem:lemma_main_1b]{3.3 (b)}).

\noindent Next, note
\[
\mathbb{E}\!\left[\frac{1}{n}\sum_{t=1}^{n} u_t^2\right]
= \frac{1}{n}\sum_{t=1}^{n} 
\exp\!\left( \frac{(1-\phi_n^{2t})\alpha^2}{2(1-\phi_n^2)} \right) = m_n,
\]
and for the second moment,
\begin{align*}
\frac{1}{n^2}\,
\mathbb{E}\!\left[\left(\sum_{t=1}^{n} u_t^2\right)^2\right]
&= \frac{\kappa_{\varepsilon}}{n^2}\sum_{t=1}^{n} \mathbb{E}[\sigma_t^4]
\;+\; \frac{2}{n^2}\sum_{1\le s<t\le n} \mathbb{E}[\sigma_s^2 \sigma_t^2],\\
&\le \frac{\kappa_{\varepsilon}}{n}\,O\left(r_n^{\frac{\alpha^2}{d}}\right)
\;+\;  e^{2\alpha^2 A_n}(1+o(1))\qquad \text{from Lemma}~\ref{lemma:main_lemma_1} \,\, \text{and}  \,\, \mathbb{E}[\varepsilon_t^4]=\kappa_{\varepsilon}.\\
& = o(1) +  e^{2\alpha^2 A_n} (1+o(1))
\end{align*}
since $\mathbb{E}[\varepsilon_t^2]=1$ and $\mathbb{E}[\varepsilon_t^4]=\kappa_{\varepsilon} < \infty$ for the noise $\{\varepsilon_t\}$, as assumed.
Since we are interested in studying the asymptotic behavior of $\sum y_t^2$, it is equally interesting to understand the results of the large sample on $\sum u_t^2$. Now it follows,
\[
\operatorname{Var}\!\left(
\frac{1}{n m_n}\sum_{t=1}^{n } u_t^2
\right)
\leq \frac{ \exp\!\big( 2\,\alpha^2\,A_n \big) (1+o(1)) - m_n^2 + o(1)}{m_n^2}
\]

\[
= \frac{e^{2 A_n \alpha^2}(1+o(1))}{e^{2 A_n \alpha^2}\big ((1+o(1)\big)^2}
\;-\;1
\;\longrightarrow\; 0.
\]
Hence, it follows that the average volatility of the process stabilizes after this normalization
\[
\left(\frac{1}{nm_n} \sum_{t=1}^{n} u_t^2\right)   \xrightarrow{p} 1. \\[10pt]
\]
Adding this into Equation~\eqref{eq:equ5_main}, this implies that 
\begin{align}\label{eq:equ6}
\frac{1}{n k_n  m_n} \sum_{t=1}^{n} y_{t-1}^2 
\xrightarrow{p} \frac{1}{2c}, 
\end{align}
because 
$1 - \rho_n^2 
= \frac{c}{k_n}\left( 2 - \frac{c}{k_n} \right) 
= \frac{2c}{k_n}\left( 1 + o(1) \right)$.

We now derive the limiting distribution of the sample covariance term
\(\sum_{t=1}^n y_{t-1}u_t\). For this, we define the martingale difference matrix and the corresponding $\sigma$ algebra, where $\alpha$ (the volatility dispersion parameter) is as defined in Section~\eqref{section-methods}. This is the main martingale step in the proof. Define the natural filtration
\(
\mathcal{F}_t
:=
\sigma\big(y_0,\varepsilon_1,\ldots,\varepsilon_t,\eta_1,\ldots,\eta_t\big),
\)
and set
\[
\tau_{nt}
:=
e^{-\alpha^2/4}
\frac{1}{\sqrt{n k_n}\,m_n}\,y_{t-1}u_t,
\qquad 1\le t\le n.
\]
Since \(\varepsilon_t\) is centered and independent of \(\mathcal{F}_{t-1}\), 
\(\{\tau_{nt},\mathcal{F}_t\}\) is a martingale difference array. Its predictable quadratic variation satisfies
\begin{align*}
\sum_{t=1}^{n}
\mathbb{E}\!\left[\tau_{nt}^2 \mid \mathcal{F}_{t-1}\right]
&=
\frac{e^{-\alpha^2/2}}{n k_n m_n^2}
\sum_{t=1}^{n}
\mathbb{E}\!\left[y_{t-1}^2 u_t^2 \mid \mathcal{F}_{t-1}\right]
\\[4pt]
&=
\frac{e^{-\alpha^2/2}}{n k_n m_n^2}
\sum_{t=1}^{n}
y_{t-1}^2\, \mathbb{E}\!\left[\sigma_t^2 \mid z_{t-1}\right]
\quad\text{(since $u_t=\sigma_t\varepsilon_t$, $\varepsilon_t\!\perp\!\mathcal{F}_{t-1}$)}
\\[4pt]
&=
\frac{1}{n k_n m_n^2}
\sum_{t=1}^{n}
y_{t-1}^2 \,\sigma_{t-1}^{2\phi_n}
\qquad( \mathbb{E}[\sigma_t^2\mid z_{t-1}]
=
\mathbb{E}\!\left[e^{\phi_n z_{t-1}+\eta_t}\mid z_{t-1}\right]
=
e^{\alpha^2/2}\sigma_{t-1}^{2\phi_n}).
\end{align*}

\noindent 
The following lemma is the key stochastic-volatility step.

\begin{lemma}\label{lem:weighted_quad_lln_main}
For the model in \eqref{eq:first model} under Assumptions~\hyperref[assm:A1]{A1} and~\hyperref[assm:A3]{A3},
\begin{align}\label{eq:eq14}
\frac{1}{n k_n m_n^2}
\sum_{t=1}^{n}y_t^2\sigma_{t-1}^{2\phi_n}
\xrightarrow{p}
\frac{e^{-\alpha^2/2}}{2c}.
\end{align}
\end{lemma}

The role of Lemma~\ref{lem:weighted_quad_lln_main} is to replace the ordinary denominator by the correct volatility-weighted quadratic form. This is where the stochastic-volatility environment differs from the homoskedastic proof. The conditional variance of the martingale score depends on the lagged volatility path, so one must control the joint behavior of \(y_t\) and \(\sigma_{t-1}^{2\phi_n}\). When $\phi_n \to 1$ slowly, the near-nonstationary volatility induces a multiplicative dispersion correction in the normalization, yielding the attenuation factor $e^{-\alpha^2/2}$. Crucially, the same factor appears in the corresponding normalization via $\tau_{nt}$, and therefore cancels. This cancellation is the mechanism behind the volatility-robust Gaussian limit. Since the difference between using \(y_t^2\) and \(y_{t-1}^2\) is asymptotically negligible under the preceding denominator bounds, Lemma~\ref{lem:weighted_quad_lln_main} implies that
\[
\sum_{t=1}^{n}
\mathbb{E}\!\left[\tau_{nt}^2\mid\mathcal{F}_{t-1}\right]
\xrightarrow{p}
\frac{e^{-\alpha^2/2}}{2c}.
\]
Next, the Lindeberg condition for the martingale array holds. For every \(\eta>0\),
\[
\sum_{t=1}^{n}
\mathbb{E}\!\left[
\tau_{nt}^2\mathbf{1}\{|\tau_{nt}|>\eta\}
\mid \mathcal{F}_{t-1}
\right]
\xrightarrow{p}0.
\]
This follows from the finite fourth moment of \(\varepsilon_t\), the initialization condition \(\mathbb{E}[y_0^4]=o(k_n^2)\), the lognormal moment bounds for \(\sigma_t\), and the rate restrictions in Assumption~\hyperref[assm:A3]{A3}. The details can be found in the \hyperref[sec:appendix]{Appendix}. Therefore, by the martingale central limit theorem,
\[
\sum_{t=1}^{n}\tau_{nt}
\xRightarrow{d}
N\!\left(0,\frac{e^{-\alpha^2/2}}{2c}\right).
\]
Equivalently,
\begin{align}\label{eq:rho_less_1_final_result}
\frac{1}{\sqrt{n k_n}\,m_n}
\sum_{t=1}^{n}y_{t-1}u_t
\xRightarrow{d}
N\!\left(0,\frac{1}{2c}\right).
\end{align}

\begin{theorem}\label{thm:thm1}
For the model in \eqref{eq:first model} under the Assumptions~\hyperref[assm:A1]{A1} and~\hyperref[assm:A3]{A3}, the following large-sample results hold:
\begin{enumerate}
\renewcommand{\labelenumi}{(\alph{enumi})}
\item 
$\displaystyle 
\frac{\sum_{t=1}^{n} y_{t-1}^2}{n\,m_n\,k_n}
\;\xrightarrow{p}\;
\frac{1}{2c},
$

\item 
$\displaystyle
\frac{1}{\sqrt{n k_n}\,m_n}
\sum_{t=1}^n y_{t-1}u_t
\xRightarrow{d}
N\!\left(0,\frac{1}{2c}\right),
$

\item 
$\displaystyle
\sqrt{n k_n}\bigl(\widehat{\rho}_n - \rho_n\bigr)
\xRightarrow{d}
N(0,2c).
$
\end{enumerate}
\end{theorem}
Theorem~\ref{thm:thm1} shows that the nearly stationary moderate-deviation limit is preserved under nearly nonstationary stochastic volatility. The condition \(k_n=o(n)\) is the same as in the homoskedastic theory, while Assumption~\hyperref[assm:A3]{A3} adds the rate separation needed to prevent volatility persistence from dominating the mean-root signal. The resulting Gaussian limit is free of the volatility parameters, so inference for \(\rho_n\) can be conducted without separately estimating the nuisance volatility scale $\alpha$, making the procedure invariant to the
volatility-dispersion specification.

We must note that no Gaussianity is required for the innovation \(\varepsilon_t\). It is enough to impose the standard moment conditions with \(\mathbb{E}[\varepsilon_t^2]=1\), and \(\mathbb{E}[\varepsilon_t^4]<\infty\), for the moderate-deviation arguments. The  moderate-deviation limits are also insensitive to the initial condition under the initialization \(\mathbb{E}[y_0^4]=o(k_n^2)\), in contrast to genuinely explosive AR(1) models where the initial value can affect the limiting distribution \citep{anderson1959asymptotic}.

Theorem~\ref{thm:thm1} gives the usual moderate-deviation interpolation of rates. If \(k_n=n^a\), \(a\in(0,1)\), then
\[
\sqrt{n k_n}\,(\widehat{\rho}_n-\rho_n)
=
n^{(1+a)/2}(\widehat{\rho}_n-\rho_n)
\xRightarrow{d} N(0,2c),
\]
which ranges from the stationary \(\sqrt n\) rate as \(a\to0\) to the local-to-unity \(n\) rate as \(a\to1\). We can also recover the corresponding stationary benchmark. For a fixed homoskedastic AR(1) with \(\rho<1\), the usual limit theory gives
\(
\sqrt n(\widehat\rho-\rho)
\xRightarrow{d}
N(0,1-\rho^2).
\)
If we formally substitute the moderate-deviation root \(\rho_n=1-c/k_n\), then
\(
1-\rho_n^2
=
\frac{2c}{k_n}-\frac{c^2}{k_n^2}
=
\frac{2c}{k_n}\{1+o(1)\}.
\)
Thus the stationary normalization \(\sqrt n/(1-\rho_n^2)^{1/2}\) becomes asymptotically proportional to \(\sqrt{n k_n}\), matching the rate in Theorem~\ref{thm:thm1}.

Our method also covers moderately negative roots. If
\(\rho_n=-1+c/k_n\), define \(y_t^*=(-1)^t y_t\) and
\(u_t^*=(-1)^t u_t\). Then
\[
\sum_{t=1}^n (y_t^*)^2=\sum_{t=1}^n y_t^2, \qquad \mathbb{E}\!\left[(u_t^*)^2\mid \mathcal F_{t-1}\right]
=
e^{\alpha^2/2}\sigma_{t-1}^{2\phi_n} \qquad \sum_{t=1}^n y_{t-1}u_t
=
-\sum_{t=1}^n y_{t-1}^*u_t^*.
\]

Thus, the transformed process has near-stationary mean dynamics with the same stochastic-volatility structure, invariant quadratic forms, while the score changes only by sign. Hence, the large-sample results and Gaussian limits carry over verbatim, up to this trivial sign change. This shows that precision is governed by the distance of the root from the unit circle, not by its sign, as in the homoskedastic setting.

The limit theory in Theorem~\ref{thm:thm1} suggests an immediate confidence interval for the autoregressive coefficient. In the nearly stationary regime,
the asymptotic \(100\lambda\%\) confidence interval for \(\rho_n\) is
\begin{equation}\label{eq:rho-infeasible-ci}
\widehat{\rho}_n
\;\pm\;
z_\lambda
\left(\frac{2c}{n k_n}\right)^{1/2},
\qquad
z_\lambda
:=
\Phi^{-1}\!\left(\frac{1+\lambda}{2}\right),
\end{equation}
where \(\Phi^{-1}\) is the standard normal quantile function. This interval is infeasible because \(c\) and \(k_n\) are unknown. However, the parametric structure in Assumption~\hyperref[assm:A1]{A1} permits one to apply likelihood-type arguments to estimate the local deviation rate from the data. The local deviation rate measures how quickly the autoregressive root approaches unity; equivalently, it is the exponent \(\gamma_n\) in the approximation \(|\rho_n-1|\approx n^{-\gamma_n}\). Combining the arguments in \citet{phillips2023estimation}, we obtain feasible confidence intervals based on the least-squares estimator \(\widehat\rho_n\) and the induced local-deviation estimator \(\widehat\gamma_n\).

\begin{theorem}
\label{thm:rho_gamma_ci}
Suppose the autoregressive root follows the moderate-deviation parameterization as in Assumption~\hyperref[assm:A1]{A1}. Let \(\widehat\rho_n\) denote the least-squares estimator of \(\rho_n\), and define
\begin{equation}
\widehat\gamma_n
:=
-\frac{\log|\widehat\rho_n-1|}{\log n}.
\label{eq:rho-gamma-hat}
\end{equation}
Let \(\lambda\in(0,1)\) denote the desired confidence level. If \(\widehat\rho_n<1\), an asymptotic \(100\lambda\%\) confidence interval for \(\rho_n\) is
\begin{equation}
\widehat\rho_n
\;\pm\;
z_\lambda\,
\frac{\sqrt{2}}{n^{(1+\widehat\gamma_n)/2}}.
\label{eq:CI-rho-MIR}
\end{equation}
The corresponding asymptotic \(100\lambda\%\) confidence interval for the local deviation rate \(\gamma_n\) is
\begin{equation}
\widehat\gamma_n
\;\pm\;
z_\lambda\,
\frac{\sqrt{2}}{n^{(1-\widehat\gamma_n)/2}\log n}.
\label{eq:CI-gamma-MIR}
\end{equation}
\end{theorem}
The quantity \(\widehat\gamma_n\) indexes the estimated local deviation of \(\widehat\rho_n\) from unity and therefore summarizes the speed at which the root approaches the unit-root boundary from below. From the inference perspective, \(\widehat\rho_n\) and \(\widehat\gamma_n\) play complementary roles. The estimator \(\widehat\rho_n\) locates the autoregressive root relative to unity, while \(\widehat\gamma_n\) measures the estimated rate of convergence toward unity. If the confidence interval for \(\rho_n\) lies strictly below one, the data provide asymptotic evidence of a non-explosive autoregressive root. Conditional on that, a confidence interval for \(\gamma_n\) contained entirely below 1 supports a mildly integrated root more strongly in the moderate-deviation near-stationary regime. Lower values of \(\widehat\gamma_n\) indicate stronger evidence that the root is away from unity, whereas values closer to one indicate weaker evidence and greater persistence near the unit-root boundary.


\section{Double local-to-unity in mildly explosive regime}\label{sec:section4}

In this section, we consider the model~\eqref{eq:first model} under Assumptions~\hyperref[assm:A2]{A2}--\hyperref[assm:A3]{A3}, so that
\[
\rho_n = 1+\frac{c}{k_n}, 
\qquad c>0.
\]
This is the mildly explosive regime. Surprisingly, the mildly explosive moderate-deviation limit requires only \(y_0=o_p(k_n^{1/2})\), and no higher moment conditions of $y_0$. In contrast to the nearly stationary case, the normalized quadratic variation \(\sum_{t=1}^{n}y_t^2\) no longer converges to a deterministic constant. Instead, the explosive mean dynamics generate a terminal stochastic scale. We therefore use the long-run volatility normalization \(\ell_n\), defined in Section~\ref{sec:our_conts}. Under \(\phi_n=1-d/\log r_n\),
\(
1-\phi_n^2=\frac{2d}{\log r_n}\{1+o(1)\},\) and hence, \(
\ell_n\sim r_n^{\alpha^2/(4d)}.
\)
Thus, \(\ell_n\) plays in the mildly explosive theory the same role that \(m_n\) plays in the nearly stationary theory and this must be tracked before the self-normalized Cauchy limit can emerge. It captures the long-run scale generated by persistent uncertainty shocks; because this scale enters both the normalized covariance and quadratic variation, it cancels from the serial-correlation statistic, leaving a Cauchy limit robust to the strength of volatility persistence.
Now we present a few interesting lemmas that isolate the main technical ingredients of the proof. 
\begin{lemma}\label{lemma:lemma1_above_1_rho}
    For the model in \eqref{eq:first model} under the Assumptions~\hyperref[assm:A2]{A2}--\hyperref[assm:A3]{A3},
\begin{align}
\rho_n^{-n} \;=\; o\!\left( \frac{k_n^t}{n}\, r_n^{-\gamma} \right)
\qquad \forall\,\gamma \in \mathbb{R},\,\, \forall\, t\geq1.
\end{align}
\end{lemma}
Lemma~\ref{lemma:lemma1_above_1_rho} provides the exponential-dominance bound needed in the mildly explosive proof. It shows that the explosive factor \(\rho_n^n\) dominates the polynomial and volatility-scale terms generated by \(k_n\), \(r_n\), and the moment bounds. This is what makes boundary terms, volatility-induced remainders, and triangular double sums negligible after explosive normalization. Thus, the lemma plays the same pivotal role as the exponential bounds in \cite{phillips2007limit}, but is strengthened to accommodate the additional stochastic-volatility scales.

\begin{lemma}\label{lemma:lemma2_above_1_rho}
For the model in~\eqref{eq:first model} under Assumptions~\hyperref[assm:A2]{A2}--\hyperref[assm:A3]{A3},
\[
\frac{\rho_n^{-n}}{\ell_n k_n}
\sum_{t=1}^{n}\sum_{j=t}^{n}\rho_n^{\,t-j-1}u_j u_t
\;\xrightarrow{L^1}\;0.
\]
\end{lemma}

Lemma~\ref{lemma:lemma2_above_1_rho} removes the triangular remainder terms that arise when the mildly explosive recursion is expanded backward. This is the step that allows the leading covariance and quadratic-variation terms to be represented through weighted sums of innovations, rather than through more complicated double sums. Our objective is to show that, after the explosive and volatility normalizations, the leading components satisfy
\[
\left(
\frac{\rho_n^{-n}\sum_{t=1}^{n} y_{t-1}u_t}{\ell_n k_n},
\;
\frac{\rho_n^{-2n}\sum_{t=1}^{n} y_{t-1}^{2}}{\ell_n k_n^{2}}
\right)
\xRightarrow{d}
(WV,W^2/2c ),
\]
where \(W\) and \(V\) are independent \(\mathcal{N}(0,1/(2c))\) random variables. To isolate these leading terms, define
\begin{align}\label{eq:def_wv_rewrite}
W_n
:=
\frac{1}{\sqrt{\ell_n k_n}}
\sum_{j=1}^{n}\rho_n^{-j}u_j,
\qquad
V_n
:=
\frac{1}{\sqrt{\ell_n k_n}}
\sum_{j=1}^{n}\rho_n^{-(n-j+1)}u_j .
\end{align}
 After Lemma~\ref{lemma:lemma2_above_1_rho} eliminates the negligible double-sum remainder, the normalized covariance $\rho_n^{-n} \sum_{t=1}^n y_{t-1}u_t/(\ell_n k_n)$ factors as a product $W_n V_n$ up to a negligible remainder, and the normalized quadratic variation $\rho_n^{-2n} \sum_{t=1}^n y_{t-1}^2/(\ell_n k_n^2)$ is asymptotically \(W_n^2/2c\). The next step is to establish the joint Gaussian limit of \((W_n,V_n)\). By the Cramér--Wold device, it is sufficient to show that, for every \(a,b\in\mathbb{R}\),
\[
aW_n+bV_n
\xRightarrow{d}
\mathcal{N}\!\left(0,\frac{a^2+b^2}{2c}\right).
\]
This yields the following lemma.

\begin{lemma}\label{lemma:wn_vn_limit_rewrite_now}
Let \((W_n,V_n)\) be defined as in~\eqref{eq:def_wv_rewrite}. Then
\[
(W_n,V_n)\xRightarrow{d}(W,V),
\]
where \(W\) and \(V\) are independent \(\mathcal{N}(0,1/(2c))\) random variables.
\end{lemma}

\noindent The proofs of Lemmas~\ref{lemma:lemma2_above_1_rho} and~\ref{lemma:wn_vn_limit_rewrite_now} are deferred to the \hyperref[sec:appendix]{Appendix}. We now use these ingredients to derive the limiting distribution of the serial-correlation statistic in the mildly explosive regime.

\noindent Returning to the time series process in~\eqref{eq:first model}, summing the identity
\[
y_t^2-\rho_n^2y_{t-1}^2
=
2\rho_n y_{t-1}u_t+u_t^2
\]
over \(t=1,\ldots,n\) and suitable normalization gives,
\begin{align*}
\frac{\rho_n^{-2n}}{k_n^2\ell_n}\sum_{t=1}^{n}y_{t-1}^2
&=
\frac{1}{k_n^2\ell_n(\rho_n^2-1)}
\left[
\rho_n^{-2n}(y_n^2-y_0^2)
-2\rho_n^{-2n+1}\sum_{t=1}^{n}y_{t-1}u_t
-\rho_n^{-2n}\sum_{t=1}^{n}u_t^2
\right].
\end{align*}
The initial condition is negligible because \(y_0=o_p(\sqrt{k_n})\) implies \(y_0/\sqrt{k_n\ell_n}=o_p(1)\). We next show that, after this normalization, only the terminal term \(y_n^2\) contributes to the limit. By Lemma~\ref{lemma:lemma1_above_1_rho} and the moment bounds for \(u_t\),
\[
\frac{\rho_n^{-2n}}{k_n\ell_n}\sum_{t=1}^{n}u_t^2=o_p(1).
\]
Similarly, using the autoregressive representation of \(y_{t-1}\),
\begin{align*}
\frac{\rho_n^{-2n+1}}{k_n\ell_n}
\sum_{t=1}^{n}y_{t-1}u_t
&=
\frac{\rho_n^{-n}y_0}{\sqrt{k_n\ell_n}}
\left(
\frac{1}{\sqrt{k_n\ell_n}}
\sum_{t=1}^{n}\rho_n^{-(n-t+1)}u_t
\right) \\
&\quad+
\frac{\rho_n^{-2n+1}}{k_n\ell_n}
\sum_{t=1}^{n}
\left(
\sum_{j=1}^{t-1}\rho_n^{t-1-j}u_j
\right)u_t
=o_p(1),
\end{align*}
where the first term is negligible by the asymptotic normality of \(V_n\), and the second term is controlled by the same exponential-dominance and volatility-moment bounds. Therefore,
\begin{align}\label{eq:eq31_main}
\frac{\rho_n^{-2n}}{\ell_n k_n^2}
\sum_{t=1}^{n}y_{t-1}^2
&=
\frac{1}{k_n(\rho_n^2-1)}
\left(
\frac{\rho_n^{-n}y_n}{\sqrt{k_n\ell_n}}
\right)^2
+o_p(1) \nonumber\\
&=
\frac{1}{k_n(\rho_n^2-1)}
\left(
\frac{y_0}{\sqrt{k_n\ell_n}}
+
\frac{1}{\sqrt{k_n\ell_n}}\sum_{j=1}^{n}\rho_n^{-j}u_j
\right)^2
+o_p(1) \nonumber\\
&\xRightarrow{d}
\frac{1}{2c}W^2,
\end{align}
where \(W\sim N(0,1/2c)\). Similarly, the normalized covariance term satisfies
\begin{align*}
\frac{\rho_n^{-n}}{k_n\ell_n}
\sum_{t=1}^{n}y_{t-1}u_t
&=
\frac{\rho_n^{-n}}{k_n\ell_n}
\sum_{t=1}^{n}y_0\rho_n^{t-1}u_t
+
\frac{\rho_n^{-n}}{k_n\ell_n}
\sum_{t=1}^{n}\sum_{j=1}^{t-1}
\rho_n^{t-1-j}u_j u_t .
\end{align*}
The first term is \(o_p(1)\). For the second term,
\begin{align}\label{eq:eq32_main}
\frac{\rho_n^{-n}}{k_n\ell_n}
\sum_{t=1}^{n}\sum_{j=1}^{t-1}
\rho_n^{t-1-j}u_j u_t
&=
\frac{\rho_n^{-n}}{k_n\ell_n}
\sum_{t=1}^{n}\sum_{j=1}^{n}
\rho_n^{t-1-j}u_j u_t
-
\frac{\rho_n^{-n}}{k_n\ell_n}
\sum_{t=1}^{n}\sum_{j=t}^{n}
\rho_n^{t-1-j}u_j u_t \nonumber\\
&=
W_nV_n+o_p(1)
\xRightarrow{d}
WV,
\end{align}
where Lemma~\ref{lemma:lemma2_above_1_rho} removes the triangular remainder and
\((W_n,V_n)\Rightarrow(W,V)\), with \(W,V\) independent \(N(0,1/(2c))\). Combining~\eqref{eq:eq31_main} and~\eqref{eq:eq32_main}, we obtain
\begin{align}\label{eq:eq33}
\frac{
\dfrac{\rho_n^{-n}\sum_{t=1}^{n}y_{t-1}u_t}{\ell_n k_n}
}{
\dfrac{\rho_n^{-2n}\sum_{t=1}^{n}y_{t-1}^2}{\ell_n k_n^2}
}
\xRightarrow{d}
\frac{WV}{W^2/(2c)}
=
2c\,\frac{V}{W}.
\end{align}
Since \(V/W\) is standard Cauchy, it follows that
\begin{align}\label{eq:cauchy_me}
\left(\frac{k_n\rho_n^n}{2c}\right)
(\widehat{\rho}_n-\rho_n)
\xRightarrow{d}
\mathcal{C},
\end{align}
where \(\mathcal{C}\sim\mathcal{C}(0,1)\).

\begin{theorem}\label{thm:thm2}
For the model in~\eqref{eq:first model} under Assumptions~\hyperref[assm:A2]{A2} and~\hyperref[assm:A3]{A3}, the following large-sample properties hold:
\[
\left(
\frac{\rho_n^{-n}\sum_{t=1}^{n}y_{t-1}u_t}{\ell_n k_n},
\;
\frac{\rho_n^{-2n}\sum_{t=1}^{n}y_{t-1}^2}{\ell_n k_n^2}
\right)
\xRightarrow{d}
\left(WV,\frac{W^2}{2c}\right),
\]
and consequently,
\[
\left(\frac{k_n\rho_n^n}{2c}\right)
(\widehat{\rho}_n-\rho_n)
\xRightarrow{d}
\mathcal{C},
\]
where \(W\) and \(V\) are independent \(N(0,1/(2c))\) random variables and \(\mathcal{C}\) is standard Cauchy.
\end{theorem}

Theorem~\ref{thm:thm2} connects the mildly explosive limit to the classical explosive AR(1) theory of \citet{white1958limiting} and \citet{basawa1984asymptotic}. The OLS deviation \(\widehat\rho_n-\rho_n\) converges at rate \(k_n\rho_n^n\), or \(n^a\rho_n^n\) when \(k_n=n^a\), as in \citet{phillips2007limit}. This rate bridges the local-to-unity order as \(a\to1\) and the fixed-explosive order as \(a\to0\), while retaining the Cauchy limit.

Since \(\rho_n=1+c/k_n\),
\(
\rho_n^2-1
=
\frac{2c}{k_n}+\frac{c^2}{k_n^2}
=
\frac{2c}{k_n}\{1+o(1)\}.
\)
Therefore,
\(
(\rho_n^2-1)^{-1}\sim \frac{k_n}{2c},
\)
so the normalization
\(
\frac{k_n\rho_n^n}{2c}(\widehat\rho_n-\rho_n)
\)
used in Theorem~\ref{thm:thm2} is the moderate-deviation counterpart of the classical fixed-explosive normalization
\[
\frac{\rho^n}{\rho^2-1}(\widehat\rho-\rho),
\qquad \rho>1 \ \text{fixed}.
\]The key difference is that, despite near-unit volatility persistence, the ratio statistic retains a Cauchy limit that is invariant to the volatility parameters and does not rely on Gaussian innovations. This implies that during incipient bubble phases, explosive mean dynamics dominate the details of volatility, yielding diagnostics of explosiveness that are less sensitive to tail and volatility misspecification than in fully explosive environments. As noted earlier, the mildly explosive limit requires only \(y_0=o_p(k_n^{1/2})\), hence the initial condition does not affect the limit law, as in the nearly stationary moderate-deviation.

Similarly, moderately explosive roots near \(-1\) are handled by the sign transformation \(y_t^*=(-1)^t y_t\), which maps the model to the positive mild-explosive case without changing the volatility path. Hence the same normalization gives the Cauchy limit, showing that inference depends on distance from the unit circle rather than the sign of the root.

The feasible inference for the mildly explosive case is obtained by using the same induced local-deviation estimator \(\widehat\gamma_n\) defined in~\eqref{eq:rho-gamma-hat}, but replacing the Gaussian quantile by the Cauchy quantile implied by Theorem~\ref{thm:thm2}.

\begin{theorem}
\label{thm:rho_gamma_ci_explosive}
Suppose the autoregressive root follows the mildly explosive moderate-deviation parameterization in Assumption~\hyperref[assm:A2]{A2}. Let \(\widehat\rho_n\) denote the least-squares estimator of \(\rho_n\), and let \(\widehat\gamma_n\) be defined as in~\eqref{eq:rho-gamma-hat}. Let \(\lambda\in(0,1)\) denote the desired confidence level and define
\[
z_\lambda^{\mathcal C}
:=
\Phi_{\mathcal C}^{-1}\!\left(\frac{1+\lambda}{2}\right),
\]
where \(\Phi_{\mathcal C}^{-1}\) is the standard Cauchy quantile function. If \(\widehat\rho_n>1\), an asymptotic \(100\lambda\%\) confidence interval for \(\rho_n\) is
\begin{equation}
\widehat\rho_n
\;\pm\;
z_\lambda^{\mathcal C}
\frac{2}{n^{\widehat\gamma_n}\widehat\rho_n^{\,n}}.
\label{eq:CI-rho-MER}
\end{equation}
The corresponding asymptotic \(100\lambda\%\) confidence interval for the local deviation rate \(\gamma_n\) is
\begin{equation}
\widehat\gamma_n
\;\pm\;
z_\lambda^{\mathcal C}
\frac{2}{\bigl(1+n^{-\widehat\gamma_n}\bigr)^n\log n}.
\label{eq:CI-gamma-MER}
\end{equation}
\end{theorem}

In empirical bubble-detection applications, \(\widehat\rho_n\) identifies whether the autoregressive root is statistically above unity, while \(\widehat\gamma_n\) measures the strength of the local explosive departure. Thus, a confidence interval for \(\rho_n\) lying strictly above one provides evidence of explosive price behavior, and the associated interval for \(\gamma_n\) indicates whether the episode is mild or strongly explosive within the moderate-deviation regime. In Section~\ref{sec:simulation}, we analyze these intervals to compare the degree of price exuberance across asset classes and historical episodes, separating sustained bubble-like dynamics from spurious signals driven mainly by short-lived volatility spikes.

\section{Empirical Diagnostics}\label{sec:simulation}

 This section connects the limit theory to empirical bubble diagnostics. We proceed in four steps. First, we verify by simulation that the normalized OLS statistics are close to their Gaussian and Cauchy limits under homoskedastic and stochastic-volatility designs. Second, we use several real-data examples to justify that volatility variation is quantitatively important and homoskedastic assumption is often infeasible during the episodes of stress. Third, we apply the proposed \((\widehat\rho_n,\widehat\gamma_n)\)-based diagnostics to recent global markets and study heterogeneity in price exuberance across assets. Finally, we compare our stochastic-volatility robust procedure with a homoskedastic PWY-type rule across canonical bubble episodes and recent comparison windows. The aim is to show that the method detects recognized exuberance while delivering more selective and economically plausible classifications in high-volatility environments. Code for all empirical analyses is documented on GitHub.

First, we examine the finite-sample accuracy of the Gaussian and Cauchy approximations by simulating model~\eqref{eq:first model} under both homoskedastic and nearly nonstationary volatility. We compare the simulated distributions of
\(
T_n=\sqrt{n k_n}\,(\widehat\rho_n-\rho_n),\) with $
\rho_n$ as in~\hyperref[assm:A1]{A1},
and
\(
S_n=\frac{\rho_n^{\,n}k_n}{2c}(\widehat\rho_n-\rho_n),\) with
$\rho_n$ as in~\hyperref[assm:A2]{A2}, with their theoretical limits \(N(0,2c)\) and \(\mathcal C(0,1)\), respectively. Accuracy is assessed using Kolmogorov--Smirnov (KS) diagnostics. For each design, we run \(100\) KS experiments based on \(100\) Monte Carlo replications and report the average KS statistic and the fraction of tests with \(p\)-value above \(0.05\). We set \(n=10^4\) and \(c=2\), in the homoskedastic design, \(\sigma_t\equiv1\), hence \(u_t=\varepsilon_t\), \(\varepsilon_t\sim N(0,1)\). In the stochastic-volatility design~\eqref{eq:first model}, we use \(d=1\), \(r_n=\log n\), and \(\alpha=0.1\).

\begin{table}[h!]
\centering
\begin{tabular}{lrrrr}
\toprule
& \multicolumn{2}{c}{Homoskedastic} & \multicolumn{2}{c}{Stochastic volatility} \\
\cmidrule(lr){2-3}\cmidrule(lr){4-5}
\(k_n\) & Mean KS stat & Acceptance prop. & Mean KS stat & Acceptance prop. \\
\midrule
constant         & 0.0884 & 0.95 & 0.0881 & 0.90 \\
\(\log n\)    & 0.0855 & 0.94 & 0.0845 & 0.96 \\
\(n^{0.10}\)  & 0.1143 & 0.76 & 0.1164 & 0.77 \\
\(n^{0.25}\)  & 0.0872 & 0.97 & 0.0900 & 0.92 \\
\(n^{0.50}\)  & 0.0956 & 0.88 & 0.0917 & 0.92 \\
\(n^{0.75}\)  & 0.1233 & 0.65 & 0.1282 & 0.63 \\
\(n/\log n\)  & 0.1250 & 0.62 & 0.1289 & 0.61 \\
\bottomrule
\end{tabular}
\caption{KS diagnostics for the nearly stationary statistic \(T_n\). Higher acceptance proportions indicate closer agreement with the Gaussian limit \(N(0,2c)\).}
\label{tab:near_stationary_ks}
\end{table}

\begin{table}[h!]
\centering
\begin{tabular}{lrrrr}
\toprule
& \multicolumn{2}{c}{Homoskedastic} & \multicolumn{2}{c}{Stochastic volatility} \\
\cmidrule(lr){2-3}\cmidrule(lr){4-5}
\(k_n\) & Mean KS stat & Acceptance prop. & Mean KS stat & Acceptance prop. \\
\midrule
\(\log n\)    & 0.0873 & 0.94 & 0.0908 & 0.88 \\
\(n^{0.10}\)  & 0.1112 & 0.81 & 0.1128 & 0.76 \\
\(n^{0.25}\)  & 0.0863 & 0.97 & 0.0879 & 0.94 \\
\(n^{0.50}\)  & 0.0834 & 0.96 & 0.0857 & 0.96 \\
\(n^{0.75}\)  & 0.0824 & 0.95 & 0.0863 & 0.95 \\
\(n^{0.90}\)  & 0.1108 & 0.73 & 0.1121 & 0.76 \\
\(n/\log n\)  & 0.0841 & 0.95 & 0.0870 & 0.93 \\
\bottomrule
\end{tabular}
\caption{KS diagnostics for the mildly explosive statistic \(S_n\). Higher acceptance proportions indicate closer agreement with the standard Cauchy limit.}
\label{tab:mild_explosive_ks}
\end{table}

\begin{figure}[H]
\centering
\includegraphics[width=0.95\textwidth]{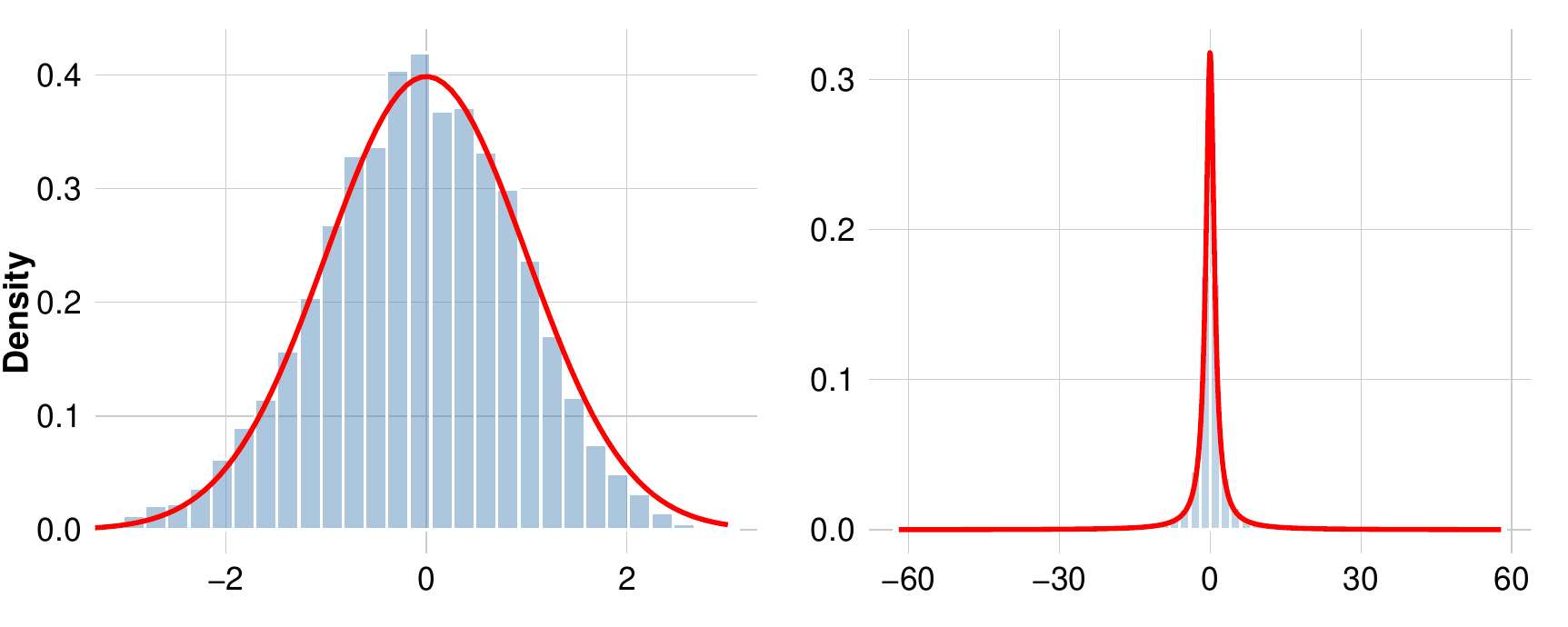}
\caption{Empirical distributions of the normalized OLS statistics under stochastic volatility. The left panel shows \(T_n=\sqrt{n k_n}(\widehat\rho_n-\rho_n)\) with the \(N(0,2c)\) overlay. The right panel shows \(S_n=\rho_n^{\,n}k_n(\widehat\rho_n-\rho_n)/(2c)\) with the standard Cauchy overlay.}
\label{fig:asymptotic_overlays}
\end{figure} 

Tables~\ref{tab:near_stationary_ks} and~\ref{tab:mild_explosive_ks} show that the simulated statistics are generally close to their predicted limits across the reported \(k_n\) regimes satisfying~\hyperref[assm:A3]{A3}. The mild deterioration for the faster \(k_n\) sequences is expected in finite samples, since these designs are closer to the unit-root boundary and require larger samples for the asymptotic approximation to dominate.

For visual diagnostics, Figure~\ref{fig:asymptotic_overlays} plots the simulated distributions of the two normalized statistics under stochastic volatility using \(B=5000\) replications. We set \(n=1000\), \(d=1\), \(\alpha=0.5\), and \(k_n=n^{0.1}\), with \(c=0.5\) in both cases. The histograms align with the limiting Gaussian and Cauchy densities, reinforcing the KS evidence that persistent stochastic volatility affects intermediate scaling but not the final self-normalized inference for \(\rho_n\).

\begin{figure}[!h]
\centering
\includegraphics[width=\linewidth]{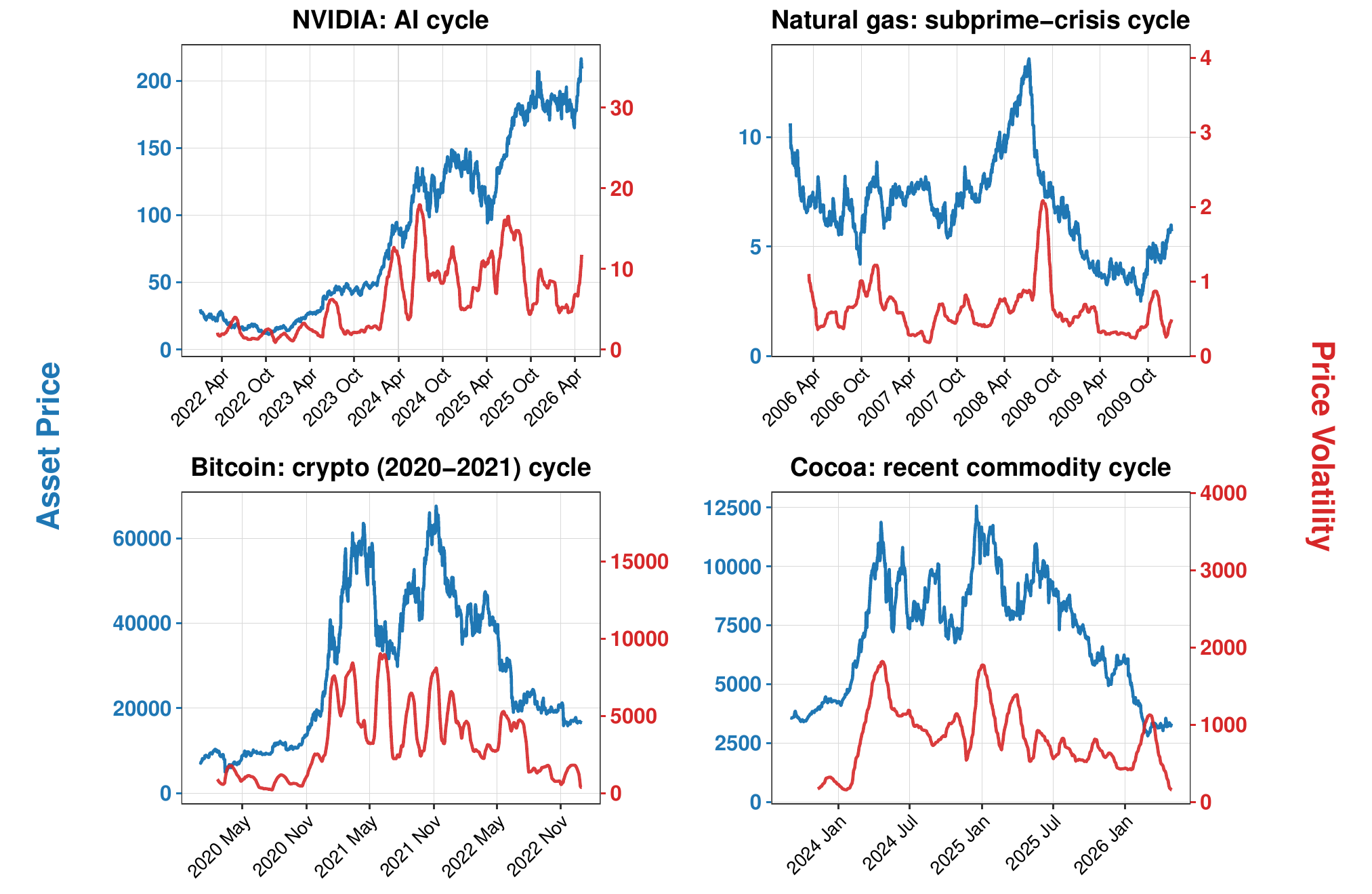}
\captionsetup{skip=2pt}
\caption{Price dynamics and rolling volatility across four episodes exhibit pronounced heteroskedasticity. Volatility is computed using a rolling eight-week window. The persistent variation in scale motivates a specification with time-varying volatility.}
\label{fig:volatility_4stocks}
\end{figure}

However, the baseline homoskedastic specification is not empirically tenable for the episodes we study now. In periods of exuberance and stress, volatility is not a background nuisance; it is part of the price dynamics itself. Large price movements are accompanied by persistent changes in uncertainty, reflecting information arrivals, speculative trading, and shifting assessments of fundamentals. Figure~\ref{fig:volatility_4stocks} illustrates this pattern for Nvidia during the recent AI cycle, natural gas around the 2007--09 crisis period, Bitcoin during the 2020--21 run-up, and cocoa during the recent commodity surge. In each case, rolling volatility varies persistently rather than remaining close to a fixed scale, motivating the time-varying volatility framework used in our theory and diagnostics.

\begin{figure}[!htbp]
    \centering
    \includegraphics[width=0.88\textwidth]{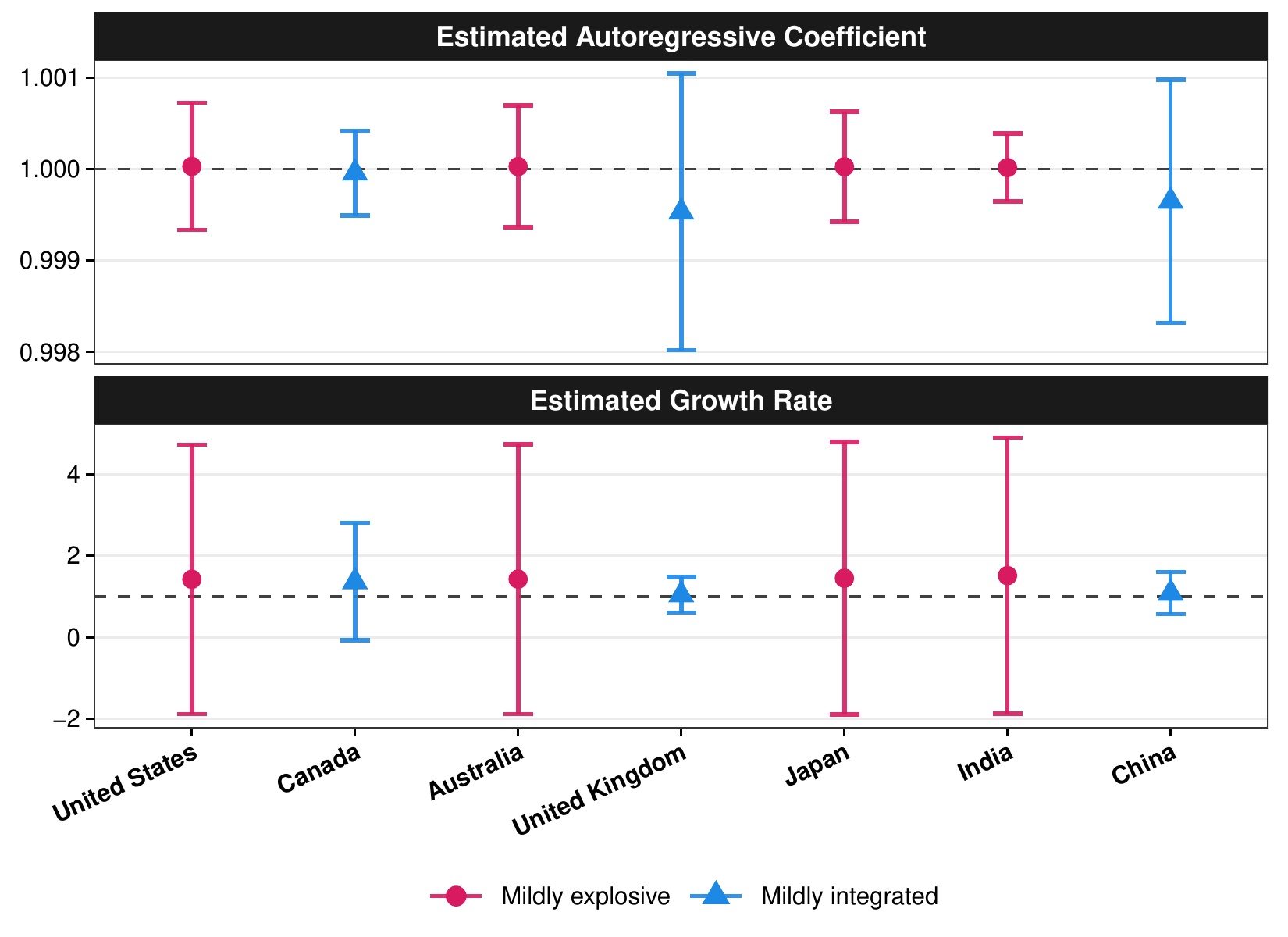}
    \caption{Estimated autoregressive coefficients and local-deviation parameters for global real-estate markets over 2020--2026. Vertical bars report asymptotic 95\% confidence intervals. The estimates are broadly clustered near the unit-root boundary, indicating high persistence across markets with limited sharply identified evidence of explosiveness.}
    \label{fig:global-re-delta-gamma}
\end{figure}

Motivated by the volatility evidence above, we apply the proposed diagnostics to recent global real-estate and commodity markets. We classify a series as mildly integrated when \(\widehat\rho_n<1\) and mildly explosive when \(\widehat\rho_n>1\).  Conditional on this classification, lower values of \(\widehat\gamma_n\) indicate a stronger local departure from unity,  while values closer to one indicate weaker evidence. For global real-estate markets over 2020--2026, the estimates are tightly clustered near unity. Australia, India, Japan, and the United States are mildly explosive with $\widehat\rho_n > 1$, see Figure~\ref{fig:global-re-delta-gamma}, but their confidence bands are wide and include one, so the evidence is weak. The United Kingdom provides the clearest non-explosive signal, with the lowest \(\widehat\gamma_n\) and the tightest confidence bound among the real-estate markets. Overall, the sector appears highly persistent, and close to the unit-root boundary.

\begin{figure}[!h]
    \centering
    \includegraphics[width=0.88\textwidth]{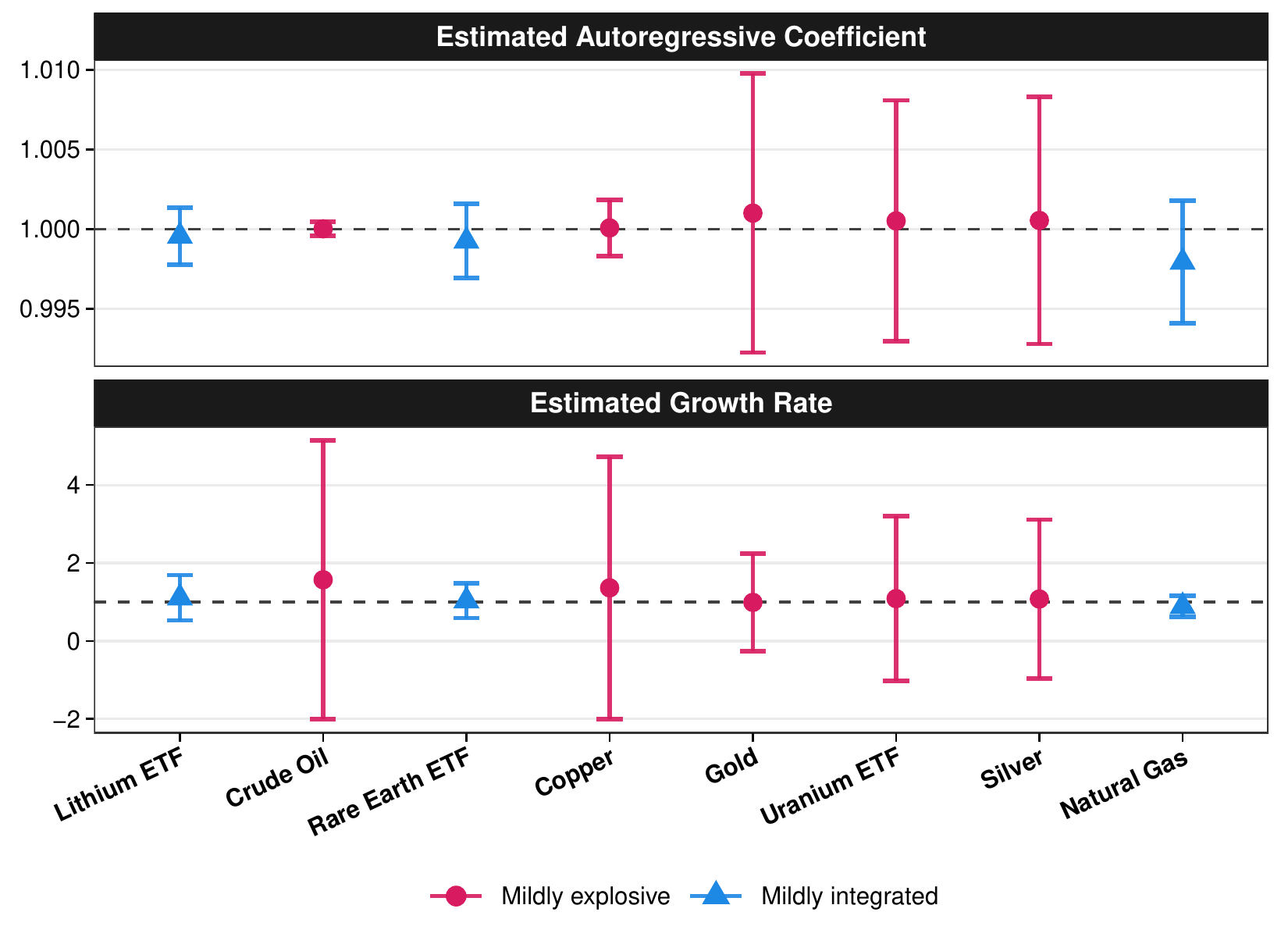}
    \caption{Estimated autoregressive coefficients and local-deviation parameters for commodity markets over 2022--2026. Vertical bars report asymptotic 95\% confidence intervals. The estimates show substantial cross-sectional heterogeneity, with stronger local-explosive evidence for some metals and weaker evidence for parts of the energy complex.}
    \label{fig:commodity-delta-gamma}
\end{figure}

For commodities over 2022--2026, our bubble-detection algorithm suggests sharper cross-sectional differences, see Figure~\ref{fig:commodity-delta-gamma}. Gold provides the strongest evidence of local explosiveness, with the largest estimated autoregressive root and a relatively low \(\widehat\gamma_n\) among the mildly explosive cases. Silver and uranium are also mildly explosive by point estimate, though their \(\gamma_n\) intervals are wider, while copper and crude oil provide weaker explosive evidence. On the non-explosive side, natural gas gives the clearest opposite signal, with \(\widehat\rho_n<1\) and the lowest \(\widehat\gamma_n=0.8867\), indicating
the strongest evidence away from unity.

\begin{figure}[H]
\centering
\includegraphics[width=0.88\linewidth]{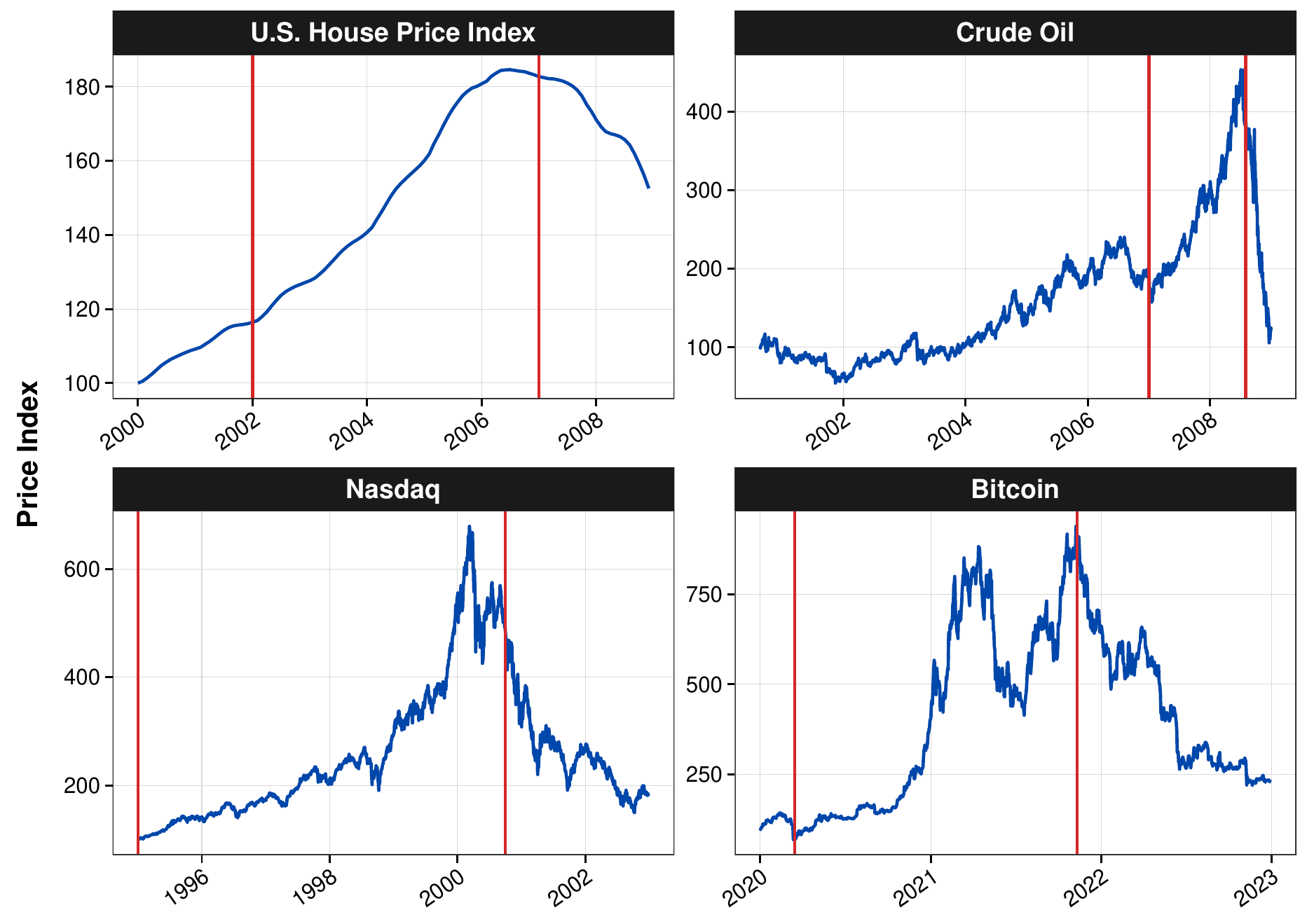}
\captionsetup{skip=0pt}
\caption{ Historical episodes of price exuberance. The panels report the U.S. house price index during the subprime housing boom, crude oil around the 2007--09 crisis period, Nasdaq during the late-1990s technology bubble, and Bitcoin during the 2020--21 crypto boom. Red vertical lines mark the detected price exuberance periods.}
\label{fig:volatility_4stocks_4_bubbles}
\end{figure}

We next apply the proposed bubble-detection algorithm to historically prominent episodes of price exuberance followed by sharp reversals. Figure~\ref{fig:volatility_4stocks_4_bubbles} reports four benchmark cases: the U.S. housing boom and crude-oil run-up around the 2007--09 subprime crisis, the late-1990s Nasdaq technology bubble, and the 2020--21 Bitcoin boom. The red vertical lines mark periods in which our method detects statistically relevant price exuberance. These examples connect to the bubble-dating framework of \cite{phillips2011dating}, while also illustrating the method in equity-technology and crypto-market episodes.  It is important to emphasize that Figure~\ref{fig:volatility_4stocks_4_bubbles} is not intended to date-stamp exact origination and collapse dates; rather, it shows that the proposed procedure detects explosive behavior during historically recognized bubble periods. For exact date-stamping of bubble origination and collapse episodes, see \cite{sarkar2026there}.

Finally, we compare the proposed stochastic-volatility robust diagnostic with the homoskedastic PWY-type rule. Although the asymptotic targets coincide with the homoskedastic benchmark, finite-sample thresholds differ once volatility specification is accounted for. For the volatility-robust rule, the cutoff is obtained by simulating the statistic under the stochastic-volatility unit-root null and using the upper \(90\%\) critical value; see GitHub for details. For the homoskedastic PWY-type rule, we use the corresponding finite-sample critical values reported in Table~1 of \citet{phillips2011explosive}. 

\begin{table}[H]
\centering
\setlength{\tabcolsep}{4pt} 
\begin{tabular}{
l
@{\hspace{1.2em}} l
@{\hspace{1.1em}} c
@{\hspace{0.8em}} c
@{\hspace{2.0em}} c
@{\hspace{0.8em}} c
}
\toprule
& & \multicolumn{2}{c}{Canonical episode}
& \multicolumn{2}{c}{Recent comparison} \\
& & & & \multicolumn{2}{c}{2020--2026} \\
\cmidrule(lr){3-4}\cmidrule(lr){5-6}
Asset & Canonical window & SV-robust & PWY & SV-robust & PWY \\
\midrule
Nikkei 225           & 1985--1989 & \BubbleYes & \BubbleYes & \BubbleNo & \BubbleYes \\
Nasdaq Composite     & 1995--2000 & \BubbleYes & \BubbleYes & \BubbleNo & \BubbleYes \\
Shanghai Composite   & 2005--2007 & \BubbleYes & \BubbleYes & \BubbleNo & \BubbleYes \\
Crude Oil            & 2007--2008 & \BubbleYes & \BubbleYes & \BubbleNo & \BubbleNo  \\
Real Estate ETF      & 2004--2007 & \BubbleYes & \BubbleYes & \BubbleNo & \BubbleYes \\
Bitcoin              & 2020--2021 & \BubbleYes & \BubbleYes & \BubbleNo & \BubbleYes \\
Ethereum             & 2020--2021 & \BubbleYes & \BubbleYes & \BubbleNo & \BubbleNo  \\
GameStop             & 2020--2021 & \BubbleNo  & \BubbleNo  & \BubbleNo & \BubbleNo  \\
ARK Innovation ETF   & 2020--2021 & \BubbleYes & \BubbleYes & \BubbleNo & \BubbleYes \\
\bottomrule
\end{tabular}
\captionsetup{skip=3pt}
\caption{Comparison of bubble classifications from the stochastic-volatility robust diagnostic and the homoskedastic PWY-type rule. Green checks indicate bubble and red crosses indicate no bubble. In the canonical windows, the two methods largely agree with historically recognized exuberance episodes. In recent windows, the PWY rule continues to flag several assets as bubbly, whereas the SV-robust diagnostic is more selective.}
\label{tab:canonical_recent_comparison}
\end{table}
In Table~\ref{tab:canonical_recent_comparison}, we first evaluate both methods on canonical or widely accepted exuberance episodes, including the Japanese Nikkei boom of 1985--1989, the Nasdaq technology bubble of 1995--2000, the Shanghai equity boom of 2005--2007, crude oil around the subprime crisis, the Bitcoin and Ethereum run-ups of 2020--2021. In these benchmark windows, our method and the PWY-type rule largely agree, both detecting the historically recognized explosive price behavior.

The recent-window comparison is more revealing. Although our method and the PWY-type rule agree on most canonical bubble episodes, they diverge sharply in recent samples. The PWY-type rule continues to classify several recent price paths as bubbles, while the volatility robust diagnostic treats them as non-explosive. This pattern suggests that a homoskedastic rule can mistake volatility movements for genuine exuberance. By accounting for time-varying volatility, the proposed diagnostic is more selective and economically plausible. It detects historically recognized bubble episodes, but cautiously avoids classifying volatility-driven movements as spurious bubbles.

\section{Conclusion}\label{sec:conclusion}

This paper develops a volatility-robust moderate-deviation theory for autoregressive models in which both the conditional mean and conditional variance can potentially approach nonstationarity. The autoregressive root satisfies \(\rho_n=1\pm c/k_n\), while the log-volatility persistence parameter also drifts toward unity. This double local-to-unity framework captures empirical environments in which prices are highly persistent and volatility is time varying and potentially unstable.

The key implication is robustness of the standard statistic. Persistent stochastic volatility can substantially amplify or compress the observed scale of price movements, especially during periods of market stress or exuberance. Nevertheless, the proposed inference for \(\rho_n\) remains valid because the procedure is built around statistics in which this volatility scale is absorbed by the normalization. As a result, the limiting inference is driven by the persistence of the autoregressive root rather than by the particular volatility-dispersion parameter. This is especially useful for empirical bubble work. Thus, the method delivers confidence intervals for persistence and local explosiveness without requiring separate estimation of the nuisance volatility parameters.

The empirical results support this robustness mechanism. Simulations show that the Gaussian and Cauchy approximations remain accurate under highly persistent stochastic volatility. Real-data applications further show that volatility is not a secondary feature during episodes of exuberance, but a part of the price dynamics itself. The proposed diagnostics detect historically recognized episodes of price exuberance, while avoiding the tendency of homoskedastic procedures to classify volatility-driven movements as spurious bubbles. This makes the method more selective and economically plausible in high-volatility environments. Methodologically, these results support self-normalized statistics and resampling methods designed for nonstationary variance environments, such as wild or multiplier bootstrap schemes in the spirit of \cite{cavaliere2008bootstrap}.

Several extensions are natural. First, the present theory complements the recursive date-stamping approach in \cite{sarkar2026there}, which can be used to study exact origination and collapse of exuberance episodes in real time. Second, an important extension is to relax correct specification of the volatility innovations and develop inference that remains reliable under weaker assumptions. Finally, multivariate and high-dimensional extensions would allow the study of common exuberance across sectors, countries, and asset classes, broadening the scope of moderate-deviation inference in settings where persistence in both prices and volatility is empirically central.

\newenvironment{acknowledgements}{%
\par\noindent\textbf{Acknowledgements}%
\begin{enumerate}%
\renewcommand\labelenumi{}%
\renewcommand\theenumi{}
}{%
\end{enumerate}%
}


\section*{AI Usage}
We used ChatGPT to assist with language editing and manuscript formatting. The scientific content is entirely the responsibility of the authors.
\newpage
\bibliographystyle{apalike}
\bibliography{bibliographyCiteDrive.bib} 

\newpage

\appendix  
\section{Appendix}\label{sec:appendix}
\numberwithin{equation}{section}
\setcounter{equation}{0}

Let us recall the time series process. We start with sequence $\{\rho_n\}$ and $\{\phi_n\}$ that determines the rate of drift of the mean and the volatility persistence coefficient to unity, respectively. Note that the autoregressive coefficient $\rho_n$ is less than 1 and slowing approaching 1. This regime is usually described as nearly stationary regime.
\subsection{Proof of Theorem~\ref{thm:thm1}}
\[
y_t=\rho_n y_{t-1}+u_t,\qquad u_t=\varepsilon_t\sigma_t,\ \ \varepsilon_t\stackrel{i.i.d.}{\sim}(0,1), \qquad \mathbb E[\varepsilon_t^4]=\kappa_\varepsilon<\infty,
\]
\[
\rho_n=1-\frac{c}{k_n},\ \ c>0,\qquad \log(\sigma_t^2)=\phi_n\log(\sigma_{t-1}^2)+\eta_t,\ \ \eta_t\stackrel{i.i.d.}{\sim}N(0,\alpha^2),
\]
\[
\phi_n=1-\frac{d}{\log r_n},\ \ d>0,\qquad 
\widehat{\rho}_n-\rho_n=\frac{\sum_{t=1}^n y_{t-1}u_t}{\sum_{t=1}^n y_{t-1}^2}.
\]

\begin{lemma}[Lemma~\ref{lemma:main_lemma_1}]\label{lemma:lemma0}
    For the model in \eqref{eq:first model} with Assumption~\hyperref[assm:A3]{A3},
    \begin{enumerate}
    \renewcommand{\labelenumi}{(\alph{enumi})}
       \item \label{e:a} $\mathbb{E}[\sigma_t^2] \leq e^{\alpha^2A_n} \;\lesssim\; r_n^{\alpha^2/(4d)},$
       \item \label{e:b} $\mathbb{E}[\sigma_t^4] \leq e^{4\alpha^2 A_n} \;\lesssim\; r_n^{\alpha^2/d},$
       \item \label{e:c} $\displaystyle \frac{1}{n^2} \sum_{s=1}^n \sum_{t \ne s}^{n} \mathbb{E}[\sigma_t^2 \sigma_s^2] \lesssim e^{2\alpha^2 A_n}.$
    \end{enumerate}
\end{lemma}

\begin{proof}

\begin{enumerate}
\renewcommand{\labelenumi}{(\alph{enumi})}

\item 
From the Gaussian AR(1) structure of \eqref{eq:vol_eq}, and with~\eqref{eq:expression_z} we have
\begin{align*}
\mathbb{E}\,[\sigma_t^{2}]
&= \exp\!\left( 
    \alpha^{2}\,\frac{1 - \phi_n^{2t}}{2\bigl(1 - \phi_n^{2}\bigr)}
  \right)
\;\le\;
\exp\!\left( 
    \alpha^{2}\,\frac{1 - \phi_n^{2n}}{2\bigl(1 - \phi_n^{2}\bigr)}
  \right) \\
&= \exp\!\left( \frac{\alpha^{2}}{2\bigl(1 - \phi_n^{2}\bigr)} \right)
   \exp\!\left( 
    - \alpha^{2}\,\frac{\phi_n^{2n}}{2\bigl(1 - \phi_n^{2}\bigr)}
  \right).
\end{align*}

Now,
\begin{align*}
\frac{\phi_n^{2n}}{1 - \phi_n^{2}}
= \frac{\bigl(1 - d / \log r_n\bigr)^{2n}}{2d / \log r_n - (d / \log r_n)^{2}}
\;\le\;
\mathcal{\mathcal{O}\!\,}\!\left(e^{-n / \log r_n}\,\log r_n\right)
= o(1).
\end{align*}

Also, note that
\begin{align*}
\exp\!\left( \frac{\alpha^{2}}{2\bigl(1 - \phi_n^{2}\bigr)} \right)
&= \exp\!\left(
    \frac{\alpha^{2}\,\log r_n}{2\bigl(2d - d^2 / \log r_n\bigr)}
  \right) \\
&= r_n^{\alpha^{2}/(4d)}\,
  \exp\!\left(
    \frac{\alpha^{2}\,\log r_n}{2}\,
    \frac{d^2 / \log r_n}{2d\,(2d - d^2 / \log r_n)}
  \right).
\end{align*}
Hence,
\begin{align*}
\exp\!\left( \frac{\alpha^{2}}{2\bigl(1 - \phi_n^{2}\bigr)} \right)
= r_n^{\alpha^{2}/(4d)} \mathcal{\mathcal{O}\!\,}\!\left(\,\exp\,\left( \frac{\alpha^2}{8}\right)\right).
\end{align*}
Hence, claim~\hyperref[e:a]{(a)} follows.

\item 
From \eqref{eq:expression_z},
\begin{align*}
z_t \;=\; \phi_n^t z_0 + \sum_{j=0}^{t-1} \phi_n^{\,j}\eta_{t-j}
&\sim N\!\left(0,\, \frac{\alpha^2(1 - \phi_n^{2t})}{1 - \phi_n^2}\right) \\
\quad\Rightarrow\quad
2z_t &\sim N\!\left(0,\, \frac{4\alpha^2(1 - \phi_n^{2t})}{1 - \phi_n^2}\right),
\end{align*}
we have the lognormal moments
\begin{align}\label{eq:equation_sigma_4_app}
\mathbb{E}[\sigma_t^4]
= \mathbb{E}\big[e^{2z_t}\big]
= \exp\!\left( 2\,\mathrm{Var}(z_t) \right)
= \exp\!\left( \frac{2\alpha^2(1 - \phi_n^{2t})}{1 - \phi_n^2} \right),
\end{align}
This implies
\begin{align*}
\mathbb{E}[\sigma_t^4]
=  e^{4\alpha^2 A_t}
\leq e^{4\alpha^2 A_n}
\;\lesssim\; r_n^{\alpha^2/d}.
\end{align*}
Hence, \hyperref[e:b]{(b)} follows.

\item 
We now measure the temporal dependence through the cross-term moments.  
Note that for $t\ge s$,
\begin{align*}
\mathbb{E}[\sigma_s^2 \sigma_t^2]
&= \mathbb{E}\big[e^{z_s+z_t}\big] \\
&= \exp\!\left\{ \tfrac{1}{2}\mathrm{Var}(z_s)
                 + \tfrac{1}{2}\mathrm{Var}(z_t)
                 + \mathrm{Cov}(z_s,z_t)\right\},
\end{align*}
with
\begin{align*}
\mathrm{Cov}(z_s,z_t)
= \alpha^2\,\phi_n^{\,t-s}\,\frac{1-\phi_n^{2s}}{1-\phi_n^2}\qquad (t\ge s),
\end{align*}
which is generally positive when $\phi_n\in(0,1)$ and $\alpha>0$.
Then the cross moment admits the representation
\begin{align*}
\mathbb{E}[\sigma_s^2 \sigma_t^2]
\;=\;
\exp\!\Big\{ \tfrac{\alpha^2}{2}\big(2A_s + 2A_t\big)
\;+\; 2\alpha^2\,\phi_n^{\,t-s}\,A_s \Big\}.
\end{align*}
Using $A_s \le A_n$, $A_t \le A_n$, and $\phi_n^{\,t-s} \le 1$, we obtain the uniform envelope
\begin{align}\label{eq:equation7_app}
\mathbb{E}[\sigma_s^2 \sigma_t^2]
\;\le\;
\exp\!\big\{ 4\,\alpha^2\,A_n \big\}.
\end{align}

However, this uniform bound is reasonably weak to guarantee the convergence of $\frac{1}{n}\sum_{t=1}^n u_t^2$ under near-nonstationary volatility, hence we require a tight bound. To obtain this, we set
\begin{equation}
\label{eq:definition_mn_app}
M_n
:=
\left\lceil
\frac{\log r_n}{d}
\log\!\Big( \frac{\log r_n}{\delta_n} \Big)
\right\rceil,
\qquad
\delta_n\downarrow0,
\qquad
\delta_n=O\!\left(\frac{1}{\log\log r_n}\right).
\end{equation}
Then, 
$$\phi_n^{M_n}
=
\left(1-\frac{d}{\log r_n}\right)^{M_n}
\le
\exp\!\left\{
-\log\!\left(\frac{\log r_n}{\delta_n}\right)
\right\}
=
\frac{\delta_n}{\log r_n}.$$
Hence, for all $h >M_n$,
\begin{align}\label{eq:prelim_bound1_app} \phi_n^h A_n
=
\mathcal{O}\,(\delta_n)
=
o(1).
\end{align}
Now note that,
\begin{align*}
\frac{1}{n^{2}}
& \sum_{1 \le s < t \le n} \mathbb{E}[\sigma_s^2 \sigma_t^2] \\
&=\frac{1}{n^{2}}
\sum_{1 \le s < t \le n}
\exp\!\left\{
\alpha^{2}(A_{s}+A_{t})
+
2 \,\alpha^{2}\phi_{n}^{\,|s-t|}A_{s}
\right\}\\
&=
\frac{1}{n^{2}}
\sum_{t=1}^{n}
\sum_{h=1}^{t-1}
\exp\!\left\{
\alpha^{2}(A_{t}+A_{t-h})
+
2\alpha^{2}\phi_{n}^{\,h}A_{t-h}
\right\}
\\[6pt]
&\le
\frac{1}{n^{2}}
\sum_{h=1}^{n}
(n-h)\,
\exp\!\left\{
2\alpha^{2}(1+\phi_{n}^{\,h})A_{n}
\right\}
\\[6pt]
&\le
\frac{1}{n}
\sum_{h=1}^{n}
\exp\!\left\{
2\alpha^{2}(1+\phi_{n}^{\,h})A_{n}
\right\}
\\[6pt]
&=
\frac{1}{n}
\left(
\sum_{h < M_{n}} e^{4\alpha^{2}A_{n}}
+
\sum_{h \ge M_{n}}
\exp\!\left\{
2\alpha^{2}(1+\phi_{n}^{\,h})A_{n}
\right\}
\right)
\\[6pt]
&\le
\frac{M_{n}}{n}e^{4\alpha^{2}A_{n}}
+
\frac{n-M_{n}}{n}
\exp\!\left\{
2\alpha^{2}(1+\mathcal{O}\!\,(\delta_n / \log r_n))A_{n} \right\}
\qquad \text{from}~\eqref{eq:prelim_bound1_app}
\\[6pt]
&\lesssim
o(1)
+
\exp\!\left\{
2\alpha^{2}A_{n}\right\} \,\, \exp\! \,(\delta_n)\qquad \text{since} \,\,M_n = o(n),\\[6pt]
& = o(1)
\;+\;
\exp\!\left\{
2\alpha^2  A_n\,\right\} (1+o(1)).
\end{align*}
Hence, the result \hyperref[e:c]{(c)} follows.
\end{enumerate}
\end{proof}

\begin{lemma}[Lemma~\ref{lemma:mean_lemma}]\label{lemma:mean_lemma_appendix_app}
For the model in \eqref{eq:first model} with~\hyperref[assm:A3]{A3}, \,\,$ m_n = e^{\alpha^2 A_n} (1+o(1)).$
\end{lemma}
\begin{proof}
From \eqref{eq:definition_mn_app}, recall that
\[
M_n
=
\left\lceil
\frac{\log r_n}{d}
\log\!\Big( \frac{\log r_n}{\delta_n} \Big)
\right\rceil,
\qquad
\delta_n\downarrow0,
\qquad
\delta_n=O\!\left(\frac{1}{\log\log r_n}\right).
\]
Since \(\phi_n=1-d/\log r_n\), we have
\begin{align}\label{eq:prelim_bound_mean_app}
\phi_n^{M_n}
&=
\left(1-\frac{d}{\log r_n}\right)^{M_n}
\le
\exp\!\left\{
-\log\!\left(\frac{\log r_n}{\delta_n}\right)
\right\}
=
\frac{\delta_n}{\log r_n}.
\end{align}
Hence, for all \(t\ge M_n\),
\[
\phi_n^t
\le
\phi_n^{M_n}
\le
\frac{\delta_n}{\log r_n}.
\]
Recall
\[
A_t
=
\frac{1-\phi_n^{2t}}{2(1-\phi_n^2)}.
\]
Then, for \(t\ge M_n\),
\[
0
\le
A_n-A_t
=
\frac{\phi_n^{2t}-\phi_n^{2n}}
     {2(1-\phi_n^2)}
\le
\frac{\phi_n^{2M_n}}
     {2(1-\phi_n^2)}.
\]
Using \eqref{eq:prelim_bound_mean_app} and
\[
1-\phi_n^2
=
\frac{2d}{\log r_n}\{1+o(1)\},
\]
we get
\[
0
\le
A_n-A_t
\le
\frac{1}{2(1-\phi_n^2)}
\left(\frac{\delta_n}{\log r_n}\right)^2
=
O\!\left(\frac{\delta_n^2}{\log r_n}\right)
=
o(1).
\]
Therefore, uniformly for \(t\ge M_n\),
$
A_t=A_n-o(1).$
Define
\(
x_t:=e^{\alpha^2 A_t}.
\)
Then, uniformly for \(t\ge M_n\),
\[
x_t
=
e^{\alpha^2 A_t}
\ge
e^{\alpha^2 A_n}\exp\{-o(1)\}
=
e^{\alpha^2 A_n}\{1+o(1)\}.
\]
Now consider
\(
S_n=\sum_{i=1}^n x_i.
\)
Since \(M_n=o(n)\), the lower bound gives
\[
S_n
\ge
\sum_{i\ge M_n}x_i
\ge
(n-M_n)e^{\alpha^2 A_n}\{1+o(1)\},
\]
and therefore
\[
m_n
=
\frac{S_n}{n}
\ge
e^{\alpha^2 A_n}\{1+o(1)\}.
\]
For the upper bound, since \(A_i\le A_n\) for all \(i\le n\),
\[
S_n
=
\sum_{i=1}^n e^{\alpha^2 A_i}
\le
n e^{\alpha^2 A_n}.
\]
Hence,
\(
m_n\le e^{\alpha^2 A_n}.
\)
Combining the lower and upper bounds yields
\begin{align}\label{eq:mu_app}
m_n
=
e^{\alpha^2 A_n}\{1+o(1)\}.
\end{align}
Finally, since
\(
A_n
=
\frac{1-\phi_n^{2n}}{2(1-\phi_n^2)}
=
\frac{1}{2(1-\phi_n^2)}+o(1),
\)
we also have
\(
m_n
=
\exp\!\left\{
\frac{\alpha^2}{2(1-\phi_n^2)}
\right\}
\{1+o(1)\}.
\)
Therefore,
\[
\frac{\alpha^2}{2(1-\phi_n^2)}
=
\frac{\alpha^2 \log r_n}{4d}
\left(1-\frac{d}{2 \log r_n}\right)^{-1}
=
\frac{\alpha^2 \log r_n}{4d}
+
\frac{\alpha^2}{8}
+
o(1).
\]
Hence,
\[
m_n
=
r_n^{\alpha^2/(4d)}
\exp\!\left\{\frac{\alpha^2}{8}+o(1)\right\}.
\] This proves one additional helpful result.
\end{proof}

\begin{lemma}[Lemma~\ref{lem:main_lemma_2} part (a)]\label{lemma:lemma1_app}
Let $y_0 = o_p(k_n^{1/2})$ and assume $r_n^{\alpha^2/(4d)} k_n = o(n)$. Then, $\frac{y_n^2}{n} \xrightarrow{p} 0$.
\end{lemma}

\begin{proof}
Write
\[
y_n = \rho_n^n y_0 + \sum_{j=1}^n \rho_n^{n-j} u_j,
\qquad u_t = \varepsilon_t \sigma_t.
\]
Then, by the elementary inequality $(a+b)^2\le 2a^2+2b^2$,
\[
\frac{y_n^2}{n}
\le \frac{2}{n}\rho_n^{2n}y_0^2 + \frac{2}{n} \Big( \sum_{j=1}^n \rho_n^{n-j} u_j \Big)^2
\le \frac{2}{n}y_0^2 + \frac{2}{n} \Big( \sum_{j=1}^n \rho_n^{n-j} u_j \Big)^2,
\]
since $\rho_n^{2n}\le 1$ for $\rho_n<1$. Using $\mathbb{E}[u_t]=0$ and independence across $t$,
\[
\frac{1}{n}\mathbb{E}\Big( \sum_{j=1}^n \rho_n^{n-j} u_j \Big)^2
= \frac{1}{n} \sum_{j=1}^n \rho_n^{2(n-j)}\mathbb{E}[\sigma_j^2].
\]
Note $z_t = \log \sigma_t^2$, $z_t = \phi_n z_{t-1} + \eta_t$. Iterating,
\[
z_t = \phi_n^t z_0 + \sum_{j=0}^{t-1}\phi_n^{\,j} \eta_{t-j}
\quad\Rightarrow\quad
\sum_{j=0}^{t-1}\phi_n^{\,j} \eta_{t-j}\sim N\!\Big(0,\frac{1-\phi_n^{2t}}{1-\phi_n^2}\alpha^2\Big).
\]
Hence,
\[
\mathbb{E}[\sigma_t^2]
= \exp\!\Big( \tfrac{1}{2}\tfrac{1-\phi_n^{2t}}{1-\phi_n^2}\alpha^2 \Big)
\le \exp\!\Big( \tfrac{\alpha^2}{2(1-\phi_n^2)} \Big)
\;\lesssim\; r_n^{\alpha^2/(4d)}, \,\, \text{from Lemma}~\hyperref[e:a]{ \text{A.1. (a)}}.\\
\]
\noindent 
Thus, we conclude that
\[
\frac{1}{n}\mathbb{E}\Big( \sum_{j=1}^n \rho_n^{n-j} u_j \Big)^2
\le \frac{1}{n} r_n^{\alpha^2/(4d)}\sum_{j=1}^n \rho_n^{2(n-j)}
= \mathcal{O}\!\,\!\Big(\frac{r_n^{\alpha^2/(4d)}k_n}{n}\Big)
= o(1),
\]
using $\sum_{j=1}^n \rho_n^{2(n-j)} \asymp (1-\rho_n^2)^{-1} = \frac{k_n}{2c}\bigl(1+o(1)\bigr)$ and $r_n^{\alpha^2/(4d)} k_n = o(n)$.
Finally, $y_0=o_p(k_n^{1/2})$ implies $\frac{y_0^2}{n}=o_p(1)$. Hence,
\[
\frac{y_n^2}{n} \xrightarrow{p} 0.
\]
\end{proof}
\begin{lemma}[Lemma~\ref{lem:main_lemma_2} part (b)]\label{lemma:lemma2_app}
If $r_n^{\frac{\alpha^2}{d}}k_n = o(n)$, then $\displaystyle \frac{1}{n} \sum_{t=1}^{n} y_{t-1} u_t \xrightarrow{p} 0.$
\end{lemma}

\begin{proof}
We have
\begin{align*}
\frac{1}{n} \sum_{t=1}^{n} y_{t-1} u_t
&= \frac{1}{n} \sum_{t=1}^{n} \left( \rho_n^{t-1} y_0 + \sum_{j=1}^{t-1} \rho_n^{t-1-j} u_j \right) u_t \\
&= \underbrace{\frac{1}{n} \sum_{t=1}^{n} \rho_n^{\,t-1} y_0 u_t}_{\equiv T_{1n}}
\;+\; \underbrace{\frac{1}{n} \sum_{t=1}^{n} \sum_{j=1}^{\,t-1} \rho_n^{\,t-1-j} u_j u_t}_{\equiv T_{2n}}.
\end{align*}
For $T_{1n}$, $\mathbb{E}[T_{1n}]=0$ and, using $y_0^2=o_p(k_n)$ and Cauchy–Schwarz,
\[
\mathbb{E}[T_{1n}^2 \mid y_0]
\;\le\; \frac{y_0^2}{n^2}\sum_{t=1}^n \rho_n^{2(t-1)} \mathbb{E}[\sigma_t^2]
\;\lesssim\;\frac{y_0^2}{n^2}\,\frac{k_n}{2c}\,r_n^{\alpha^2/(4d)}
\;=\; o_p(1).
\]
Hence $T_{1n}=o_p(1)$. For $T_{2n}$,
\begin{align*}
\mathbb{E}[T_{2n}^2]
&= \frac{1}{n^2} \sum_{t=1}^{n} \mathbb{E}\!\left[ \sum_{j=1}^{t-1} \rho_n^{\,2(t-1-j)} u_j^2 u_t^2 \right] \\
&= \frac{1}{n^2} \sum_{t=1}^{n} \sum_{j=1}^{t-1} 
\rho_n^{\,2(t-1-j)} \,\mathbb{E}[\sigma_j^2 \sigma_t^2]\, \\
&\le \frac{1}{n^2} \sum_{t=1}^{n}
\frac{1}{1-\rho_n^2}\,\,\mathcal{O}\!\,(r_n^{\alpha^2/d}) \,\,\,\, \text{from Cauchy Schwartz and Lemma}~\hyperref[e:b]{\text{A.1 (b)}}.\\ 
&= \mathcal{O}\!\,\!\left(\frac{r_n^{\alpha^2/d}\,k_n}{n}\right)
= o(1),
\end{align*}
by the assumed growth condition. Thus $T_{2n}=o_p(1)$, and the result follows.
\end{proof}

\begin{lemma}[Lemma~\ref{lem:weighted_quad_lln_main}]\label{lem:weighted_quad_lln}
For the model in~\eqref{eq:first model} under Assumptions~\hyperref[assm:A1]{A1} and~\hyperref[assm:A3]{A3},
\begin{align}\label{eq:weighted_quad_lln_app}
\frac{1}{n k_n m_n^2}
\sum_{t=1}^{n}
y_{t-1}^2\sigma_{t-1}^{2\phi_n}
\xrightarrow{p}
\frac{e^{-\alpha^2/2}}{2c}.
\end{align}
\end{lemma}

\begin{proof}
Let
\[
a_t:=\sigma_{t-1}^{2\phi_n}=e^{\phi_n z_{t-1}}.
\]
We prove~\eqref{eq:weighted_quad_lln_app} by showing convergence of the mean and concentration around the mean. First,
\begin{align}
\mathbb E[a_t]
&=
\mathbb E[e^{\phi_n z_{t-1}}]
=
\exp\left\{
\frac{\alpha^2\phi_n^2(1-\phi_n^{2(t-1)})}
{2(1-\phi_n^2)}
\right\}
=
\exp\{\alpha^2\phi_n^2 A_{t-1}\}.
\label{eq:mean-at-app}
\end{align}
For \(t\ge M_n\), Lemma~\ref{lemma:mean_lemma_appendix_app} gives \(A_t=A_n+o(1)\). Moreover,
\[
(1-\phi_n^2)A_n
=
\frac{1-\phi_n^{2n}}{2}
=
\frac12+o(1).
\]
Hence, uniformly for \(t\ge M_n\),
\begin{align}
\mathbb E[a_t]
=
\exp\{\alpha^2A_n-\alpha^2/2+o(1)\}
=
e^{-\alpha^2/2}m_n\{1+o(1)\}.
\label{eq:mean-at-asymp-app}
\end{align}
\noindent 
Using the AR representation
\[
y_{t-1}
=
\rho_n^{t-1}y_0
+
\sum_{j=1}^{t-1}\rho_n^{t-1-j}u_j,
\]
and conditioning on the volatility path,
\[
\mathbb E_\varepsilon[y_{t-1}^2]
=
\rho_n^{2(t-1)}y_0^2
+
\sum_{j=1}^{t-1}\rho_n^{2(t-1-j)}\sigma_j^2.
\]
Therefore,
\begin{align}
\mathbb E[y_{t-1}^2a_t]
=
\rho_n^{2(t-1)}\mathbb E[y_0^2a_t]
+
\sum_{j=1}^{t-1}
\rho_n^{2(t-1-j)}
\mathbb E[\sigma_j^2a_t].
\label{eq:weighted-mean-decomp-app}
\end{align}
The initial-condition term is negligible after normalization:
\begin{align}
\frac{1}{n k_n m_n^2}
\sum_{t=1}^{n}
\rho_n^{2(t-1)}
\mathbb E[y_0^2a_t]
&\le
\frac{C\mathbb E[y_0^2]}{n k_n m_n^2}
\sum_{t=1}^{n}\rho_n^{2(t-1)}m_n
=o(1),
\label{eq:init-weighted-negligible-app}
\end{align}
using \(\mathbb E[y_0^2]=o(k_n)\), \(\sum_t\rho_n^{2t}=O(k_n)\), and \(m_n\to\infty\).
\noindent 
For the double sum, put \(h=t-j\). Then
\[
\sum_{t=1}^{n}
\sum_{j=1}^{t-1}
\rho_n^{2(t-1-j)}
\mathbb E[\sigma_j^2a_t]
=
\sum_{j=1}^{n-1}
\sum_{h=1}^{n-j}
\rho_n^{2(h-1)}
\mathbb E[\sigma_j^2\sigma_{j+h-1}^{2\phi_n}].
\]
For \(h>M_n\), the covariance between \(z_j\) and \(z_{j+h-1}\) is negligible because
\[
\phi_n^{h-1}A_j=o(1).
\]
Thus, uniformly over the dominant range,
\[
\mathbb E[\sigma_j^2\sigma_{j+h-1}^{2\phi_n}]
=
\mathbb E[\sigma_j^2]\,
\mathbb E[\sigma_{j+h-1}^{2\phi_n}]
\{1+o(1)\}.
\]
For \(j\ge M_n\) and \(j+h-1\ge M_n\), Lemma~\ref{lemma:mean_lemma_appendix_app} and~\eqref{eq:mean-at-asymp-app} imply
\begin{align}
\mathbb E[\sigma_j^2\sigma_{j+h-1}^{2\phi_n}]
=
e^{-\alpha^2/2}m_n^2\{1+o(1)\}.
\label{eq:dominant-vol-moment-app}
\end{align}
The contribution from \(h\le M_n\) is negligible. By Cauchy--Schwarz and the lognormal moment bounds, for some finite \(q>0\),
\[
\mathbb E[\sigma_j^2\sigma_{j+h-1}^{2\phi_n}]
\le
C r_n^q m_n^2.
\]
Therefore,
\begin{align}
\frac{1}{n k_n m_n^2}
\sum_{j=1}^{n}
\sum_{h\le M_n}
\rho_n^{2(h-1)}
\mathbb E[\sigma_j^2\sigma_{j+h-1}^{2\phi_n}]
&\le
C\frac{M_n r_n^q}{k_n}
=o(1),
\label{eq:short-lag-negligible-app}
\end{align}
by Assumption~\hyperref[assm:A3]{A3}. The contribution from \(j<M_n\) is also \(o(1)\), since \(M_n=o(n)\). Hence,
\begin{align}
\frac{1}{n k_n m_n^2}
\sum_{t=1}^{n}
\mathbb E[y_{t-1}^2a_t]
&=
\frac{e^{-\alpha^2/2}\{1+o(1)\}}{n k_n}
\sum_{j=M_n}^{n}
\sum_{h>M_n}^{n-j}
\rho_n^{2(h-1)}
+o(1).
\label{eq:weighted-mean-before-limit-app}
\end{align}
Since \(M_n=o(k_n)\),
\[
\sum_{h>M_n}^{n-j}\rho_n^{2(h-1)}
=
\frac{\rho_n^{2M_n}}{1-\rho_n^2}\{1+o(1)\}
=
\frac{k_n}{2c}\{1+o(1)\}.
\]
Therefore,
\begin{align}
\frac{1}{n k_n m_n^2}
\sum_{t=1}^{n}
\mathbb E[y_{t-1}^2a_t]
\to
\frac{e^{-\alpha^2/2}}{2c}.
\label{eq:mean-weighted-limit-app}
\end{align}
\noindent 
It remains to show concentration. Let
\[
R_n
:=
\sum_{t=1}^{n}
\left\{
y_{t-1}^2a_t-\mathbb E[y_{t-1}^2a_t]
\right\}.
\]
It is enough to prove
\[
\mathbb E[R_n^2]=o(n^2k_n^2m_n^4).
\]
Using the AR expansion for \(y_{t-1}\), the finite fourth moment of \(\varepsilon_t\), and the lognormal moment bounds, there exists a finite \(q>0\) such that
\begin{align}
\mathbb E[(y_{t-1}^2a_t)^2]
\le
C k_n^2m_n^4 r_n^q.
\label{eq:diagonal-bound-app}
\end{align}
For \(s<t\), putting \(h=t-s\), the covariance between
\(y_{t-1}^2a_t\) and \(y_{s-1}^2a_s\) is controlled by the AR memory \(\rho_n^h\) and the volatility memory \(\phi_n^h\). The same moment bounds imply
\begin{align}
\left|
\operatorname{Cov}(y_{t-1}^2a_t,\ y_{s-1}^2a_s)
\right|
\le
C k_n^2m_n^4 r_n^q(\rho_n^h+\phi_n^h).
\label{eq:covariance-bound-app}
\end{align}
Consequently,
\begin{align}
\mathbb E[R_n^2]
&\le
C n k_n^2m_n^4r_n^q
+
C n k_n^2m_n^4r_n^q
\sum_{h=1}^{n-1}(\rho_n^h+\phi_n^h)
\nonumber\\
&\le
C n k_n^2m_n^4r_n^q
\left\{
1+\frac{1}{1-\rho_n}+\frac{1}{1-\phi_n}
\right\}.
\label{eq:Rn-second-moment-app}
\end{align}
Since
\[
1-\rho_n=\frac{c}{k_n},
\qquad
1-\phi_n=\frac{d}{\log r_n},
\]
we obtain
\[
\mathbb E[R_n^2]
\le
C n k_n^2m_n^4r_n^q(1+k_n+\log r_n).
\]
Thus,
\begin{align}
\frac{\mathbb E[R_n^2]}{n^2k_n^2m_n^4}
\le
C r_n^q\frac{1+k_n+\log r_n}{n}
=o(1),
\label{eq:Rn-variance-small-app}
\end{align}
by Assumption~\hyperref[assm:A3]{A3}. Hence
\[
\frac{R_n}{n k_n m_n^2}=o_p(1).
\]
Combining this with~\eqref{eq:mean-weighted-limit-app} proves~\eqref{eq:weighted_quad_lln_app}.
\end{proof}
\noindent We now have all the ingredients to prove main Theorem~\ref{thm:thm1}. We first show the denominator limit. From
\[
y_t^2-\rho_n^2y_{t-1}^2
=
u_t^2+2\rho_n y_{t-1}u_t,
\]
summing over \(t=1,\ldots,n\) and dividing by \(n\) gives
\begin{align}\label{eq:equ5_app}
(1-\rho_n^2)\frac{1}{n}\sum_{t=1}^{n}y_{t-1}^2
&=
\frac{1}{n}(y_0^2-y_n^2)
+
\frac{1}{n}\sum_{t=1}^{n}u_t^2
+
2\rho_n\frac{1}{n}\sum_{t=1}^{n}y_{t-1}u_t
\nonumber\\
&=
\frac{1}{n}\sum_{t=1}^{n}u_t^2+o_p(1),
\end{align}
where we use \(y_0=o_p(k_n^{1/2})\), Lemma~\ref{lemma:lemma1_app}, and Lemma~\ref{lemma:lemma2_app}. By the volatility moment bounds and Lemma~\ref{lemma:mean_lemma_appendix_app},
\[
\frac{1}{n m_n}\sum_{t=1}^{n}u_t^2
\xrightarrow{p}1.
\]
Since
\[
1-\rho_n^2
=
\frac{2c}{k_n}\{1+o(1)\},
\]
we obtain
\begin{align}\label{eq:denom_limit_app}
\frac{1}{n k_n m_n}
\sum_{t=1}^{n}y_{t-1}^2
\xrightarrow{p}
\frac{1}{2c}.
\end{align}
\noindent 
We next derive the limiting distribution of the martingale score. Define
\[
\tau_{nt}
=
e^{-\alpha^2/4}
\frac{1}{\sqrt{n k_n}\,m_n}
y_{t-1}u_t,
\qquad
\mathcal F_t
=
\sigma(y_0,\varepsilon_1,\ldots,\varepsilon_t,\eta_1,\ldots,\eta_t).
\]
Then \(\{\tau_{nt},\mathcal F_t\}\) is a martingale difference array. Its predictable quadratic variation satisfies
\begin{align}
\sum_{t=1}^{n}
\mathbb E[\tau_{nt}^2\mid\mathcal F_{t-1}]
&=
\frac{1}{n k_n m_n^2}
\sum_{t=1}^{n}
y_{t-1}^2\sigma_{t-1}^{2\phi_n}.
\label{eq:cond_var_tau_app}
\end{align}
By Lemma~\ref{lem:weighted_quad_lln},
\begin{align}
\sum_{t=1}^{n}
\mathbb E[\tau_{nt}^2\mid\mathcal F_{t-1}]
\xrightarrow{p}
\frac{e^{-\alpha^2/2}}{2c}.
\label{eq:tau_variance_limit_app}
\end{align}

\noindent It remains to verify the Lindeberg condition. For any \(\eta>0\), let
\[
L_n(\eta)
:=
\sum_{t=1}^{n}
\mathbb E\!\left[
\tau_{nt}^2
\mathbf 1\{|\tau_{nt}|>\eta\}
\mid\mathcal F_{t-1}
\right].
\]
Using \(x^2\mathbf 1\{|x|>\eta\}\le \eta^{-2}x^4\),
\begin{align}
L_n(\eta)
&\le
\frac{e^{-\alpha^2}}{\eta^2n^2k_n^2m_n^4}
\sum_{t=1}^{n}
y_{t-1}^4
\mathbb E[u_t^4\mid\mathcal F_{t-1}].
\label{eq:lindeberg_bound_app}
\end{align}
Since
\[
\mathbb E[u_t^4\mid\mathcal F_{t-1}]
=
\kappa_\varepsilon e^{2\alpha^2}\sigma_{t-1}^{4\phi_n},
\]
it is enough to show
\begin{align}\label{eq:weighted_fourth_target_app}
\frac{1}{n^2k_n^2m_n^4}
\sum_{t=1}^{n}
y_{t-1}^4\sigma_{t-1}^{4\phi_n}
=o_p(1).
\end{align}
The proof follows from the AR representation of \(y_{t-1}\), the finite fourth moment of \(\varepsilon_t\), the initialization condition, and the lognormal moment bounds. In particular, for some finite \(q>0\),
\[
\mathbb E[y_{t-1}^4\sigma_{t-1}^{4\phi_n}]
\le
C k_n^2m_n^4r_n^q
\]
uniformly in \(t\le n\). Hence
\[
\mathbb E\left[
\frac{1}{n^2k_n^2m_n^4}
\sum_{t=1}^{n}
y_{t-1}^4\sigma_{t-1}^{4\phi_n}
\right]
\le
C\frac{r_n^q}{n}
=o(1),
\]
by Assumption~\hyperref[assm:A3]{A3}. Thus~\eqref{eq:weighted_fourth_target_app} holds by Markov's inequality, and \(L_n(\eta)\xrightarrow{p}0\).

Combining~\eqref{eq:tau_variance_limit_app} with the Lindeberg condition, the martingale central limit theorem yields
\[
\sum_{t=1}^{n}\tau_{nt}
\xRightarrow{d}
N\!\left(0,\frac{e^{-\alpha^2/2}}{2c}\right).
\]
Since
\[
\sum_{t=1}^{n}\tau_{nt}
=
e^{-\alpha^2/4}
\frac{1}{\sqrt{n k_n}m_n}
\sum_{t=1}^{n}y_{t-1}u_t,
\]
we obtain
\begin{align}\label{eq:score_limit_app}
\frac{1}{\sqrt{n k_n}m_n}
\sum_{t=1}^{n}y_{t-1}u_t
\xRightarrow{d}
N\!\left(0,\frac{1}{2c}\right).
\end{align}

\noindent Finally, using
\[
\widehat\rho_n-\rho_n
=
\frac{\sum_{t=1}^{n}y_{t-1}u_t}
     {\sum_{t=1}^{n}y_{t-1}^2},
\]
and combining~\eqref{eq:denom_limit_app} with~\eqref{eq:score_limit_app}, we get
\begin{align}
\sqrt{n k_n}
(\widehat\rho_n-\rho_n)
\xRightarrow{d}
N(0,2c).
\end{align}
This proves Theorem~\ref{thm:thm1}.

\subsection{Proof of Theorem~\ref{thm:rho_gamma_ci}}

Consider the autoregression
\[
    y_t=\rho_n y_{t-1}+u_t,
    \qquad
    u_t=\sigma_t\varepsilon_t,
\]
under the regularity conditions of Theorem~\ref{thm:thm1}. Let
\begin{align}\label{main_rho_equation}
    \widehat\rho_n
    =
    \frac{\sum_{t=1}^n y_{t-1}y_t}
         {\sum_{t=1}^n y_{t-1}^2},
    \qquad
    \widehat\gamma_n
    =
    -\frac{\log|\widehat\rho_n-1|}{\log n}.
\end{align}
Suppose
$
    \rho_n=1-\frac{1}{n^{\gamma_n}},
    \,
    0<\gamma_n<1.
$
Then
\begin{align}\label{main_rho_equation_CI}
    n^{(1+\gamma_n)/2}
    \bigl(\widehat\rho_n-\rho_n\bigr)
    \Rightarrow Z,
    \qquad
    Z\sim N(0,2),
\end{align}
and
\begin{align}\label{main_gamma_equation_CI}
    n^{(1-\gamma_n)/2}\log n
    \bigl(\widehat\gamma_n-\gamma_n\bigr)
    \Rightarrow Z,
    \qquad
    Z\sim N(0,2).
\end{align}
Consequently, if \(\widehat\rho_n<1\), an asymptotic \(100\lambda\%\) confidence interval for \(\rho_n\) is
\[
    \widehat\rho_n
    \pm
    z_\lambda
    \frac{\sqrt{2}}{n^{(1+\widehat\gamma_n)/2}},
    \qquad
    z_\lambda
    =
    \Phi^{-1}\!\left(\frac{1+\lambda}{2}\right),
\]
and the corresponding asymptotic \(100\lambda\%\) confidence interval for \(\gamma_n\) is
\[
    \widehat\gamma_n
    \pm
    z_\lambda
    \frac{\sqrt{2}}
    {n^{(1-\widehat\gamma_n)/2}\log n}.
\]

\begin{proof}
We give the proof in the localizing-rate parametrization. Define
\begin{equation}\label{eq:near_A_defs}
    A_n:=\rho_n-1=-n^{-\gamma_n},
    \qquad
    \widehat A_n:=\widehat\rho_n-1 .
\end{equation}
By Theorem~\ref{thm:thm1},
\begin{equation}\label{eq:near_rho_clt}
    n^{(1+\gamma_n)/2}
    (\widehat\rho_n-\rho_n)
    \xRightarrow{d} Z,
    \qquad
    Z\sim N(0,2).
\end{equation}
Equivalently, for some sequence \(Z_n\xRightarrow{d}Z\),
\begin{equation}\label{eq:near_Zn_rep}
    \widehat\rho_n-\rho_n
    =
    n^{-(1+\gamma_n)/2}Z_n .
\end{equation}
Using~\eqref{eq:near_A_defs} and~\eqref{eq:near_Zn_rep},
\begin{align}
    \frac{\widehat A_n}{A_n}
    &=
    \frac{A_n+(\widehat\rho_n-\rho_n)}{A_n} \notag\\
    &=
    1
    -
    n^{\gamma_n}(\widehat\rho_n-\rho_n) \notag\\
    &=
    1
    -
    n^{-(1-\gamma_n)/2}Z_n .
    \label{eq:near_A_ratio}
\end{align}
Since \(0<\gamma_n<1\) and \(Z_n=O_p(1)\),
\begin{equation}\label{eq:near_small_term}
    n^{-(1-\gamma_n)/2}Z_n=o_p(1).
\end{equation}
Thus \(\widehat A_n/A_n=1+o_p(1)\). Moreover, from~\eqref{eq:near_A_ratio},
\begin{equation}\label{eq:near_abs_A}
    |\widehat A_n|
    =
    n^{-\gamma_n}
    \left|
    1-n^{-(1-\gamma_n)/2}Z_n
    \right|.
\end{equation}
Recalling that
\[
    \widehat\gamma_n
    =
    -\frac{\log|\widehat\rho_n-1|}{\log n}
    =
    -\frac{\log|\widehat A_n|}{\log n},
\]
and taking logarithms in~\eqref{eq:near_abs_A}, we obtain
\begin{align}
    \widehat\gamma_n-\gamma_n
    &=
    -
    \frac{
    \log\left|
    1-n^{-(1-\gamma_n)/2}Z_n
    \right|
    }{\log n}.
    \label{eq:near_gamma_exact}
\end{align}
By~\eqref{eq:near_small_term}, the argument inside the logarithm is positive with probability tending to one, and the expansion
\[
    \log(1+x)=x+o_p(x)
    \qquad (x=o_p(1))
\]
applies with \(x=-n^{-(1-\gamma_n)/2}Z_n\). Hence
\begin{align}
    \log\left|
    1-n^{-(1-\gamma_n)/2}Z_n
    \right|
    &=
    -
    n^{-(1-\gamma_n)/2}Z_n
    +
    o_p\!\left(n^{-(1-\gamma_n)/2}\right).
    \label{eq:near_log_expansion}
\end{align}
Substituting~\eqref{eq:near_log_expansion} into~\eqref{eq:near_gamma_exact} gives
\begin{equation}\label{eq:near_gamma_expansion}
    \widehat\gamma_n-\gamma_n
    =
    \frac{
    n^{-(1-\gamma_n)/2}Z_n
    }{\log n}
    +
    o_p\!\left(
    \frac{n^{-(1-\gamma_n)/2}}{\log n}
    \right).
\end{equation}
Therefore,
\begin{equation}\label{eq:near_gamma_limit}
    n^{(1-\gamma_n)/2}\log n
    (\widehat\gamma_n-\gamma_n)
    =
    Z_n+o_p(1)
    \xRightarrow{d} Z,
    \qquad
    Z\sim N(0,2).
\end{equation}
It remains to justify replacing the unknown \(\gamma_n\) by \(\widehat\gamma_n\) in the feasible normalizations. From~\eqref{eq:near_gamma_limit},
\begin{equation}\label{eq:near_gamma_rate}
    n^{(1-\gamma_n)/2}\log n
    (\widehat\gamma_n-\gamma_n)
    =
    O_p(1),
\end{equation}
so
\begin{equation}\label{eq:near_gamma_log_rate}
    (\widehat\gamma_n-\gamma_n)\log n
    =
    O_p\!\left(n^{-(1-\gamma_n)/2}\right)
    =
    o_p(1).
\end{equation}
Consequently,
\begin{align}
    \frac{n^{(1+\widehat\gamma_n)/2}}
         {n^{(1+\gamma_n)/2}}
    &=
    \exp\!\left\{
    \frac{1}{2}(\widehat\gamma_n-\gamma_n)\log n
    \right\}
    =
    1+o_p(1),
    \label{eq:near_feasible_rho_norm}
\end{align}
and similarly,
\begin{align}
    \frac{n^{(1-\widehat\gamma_n)/2}}
         {n^{(1-\gamma_n)/2}}
    &=
    \exp\!\left\{
    -\frac{1}{2}(\widehat\gamma_n-\gamma_n)\log n
    \right\}
    =
    1+o_p(1).
    \label{eq:near_feasible_gamma_norm}
\end{align}
Combining~\eqref{eq:near_rho_clt} with~\eqref{eq:near_feasible_rho_norm} yields
\[
    n^{(1+\widehat\gamma_n)/2}
    (\widehat\rho_n-\rho_n)
    \xRightarrow{d}
    N(0,2),
\]
which gives the feasible confidence interval for \(\rho_n\). Likewise, combining~\eqref{eq:near_gamma_limit} with~\eqref{eq:near_feasible_gamma_norm} gives
\[
    n^{(1-\widehat\gamma_n)/2}\log n
    (\widehat\gamma_n-\gamma_n)
    \xRightarrow{d}
    N(0,2),
\]
and hence the stated feasible confidence interval for \(\gamma_n\).
\end{proof}

\subsection{Proof of Theorem~\ref{thm:thm2}}
Let us now consider the time series process as described in Section~\ref{sec:section4}.
\[
y_t=\rho_n y_{t-1}+u_t,\qquad u_t=\varepsilon_t\sigma_t,\ \ \varepsilon_t\stackrel{i.i.d.}{\sim}(0,1), 
\]
\[
\rho_n=1+\frac{c}{k_n},\ \ c>0,\qquad \log(\sigma_t^2)=\phi_n\log(\sigma_{t-1}^2)+\eta_t,\ \ \eta_t\stackrel{i.i.d.}{\sim}N(0,\alpha^2),
\]
\[
\phi_n=1-\frac{d}{\log r_n},\ \ d>0.
\]
Note that the autoregressive coefficient $\rho_n$ is now more than 1 and slowing approaching 1. This regime is usually described as ``mild" explosive regime under nearly non-stationary stochastic volatility. We will prove the distributional results of $\widehat{\rho_n}$, and asymptotic distribution of 
\[
\widehat{\rho}_n-\rho_n=\frac{\sum_{t=1}^n y_{t-1}u_t}{\sum_{t=1}^n y_{t-1}^2}.
\]

\begin{lemma}[Lemma~\ref{lemma:lemma1_above_1_rho}]\label{lem:lemma3}
Under the model specifications described in Section~\ref{sec:section4} and Assumption~\hyperref[assm:A3]{A3},
\[
\rho_n^{-n}
=
o\!\left(
\frac{k_n^t}{n}\,r_n^{-\gamma}
\right),
\qquad
\forall\,\gamma\in\mathbb{R},\quad \forall\,t\ge 1.
\]
\end{lemma}

\begin{proof}
It is enough to show that
\[
\log\!\left(
n\,\frac{\rho_n^{-n}r_n^\gamma}{k_n^t}
\right)
\to -\infty.
\]
Since \(\rho_n=1+c/k_n\),
\[
-n\log\rho_n
=
-n\log\!\left(1+\frac{c}{k_n}\right)
=
-\frac{c n}{k_n}
+
\mathcal{O}\!\left(\frac{n}{k_n^2}\right).
\]
Therefore,
\begin{align*}
\log\!\left(
n\,\frac{\rho_n^{-n}r_n^\gamma}{k_n^t}
\right)
&=
-\frac{c n}{k_n}
+
\mathcal{O}\!\left(\frac{n}{k_n^2}\right)
+
\gamma\log r_n
+
\log n
-
t\log k_n \\
&=
\frac{c n}{k_n}
\left[
-1
+
\mathcal{O}\!\left(\frac{1}{k_n}\right)
+
\frac{k_n}{c n}\gamma\log r_n
+
\frac{k_n}{c n}\{\log n-t\log k_n\}
\right].
\end{align*}
For \(t\ge1\),
\[
\log n-t\log k_n
=
\log\!\left(\frac{n}{k_n}\right)
-(t-1)\log k_n
\le
\log\!\left(\frac{n}{k_n}\right).
\]
Since \(n/k_n\to\infty\),
\[
\frac{k_n}{n}\log\!\left(\frac{n}{k_n}\right)\to0.
\]
Also, by Assumption~\hyperref[assm:A3]{A3}, \(k_n r_n^q=o(n)\) for every fixed \(q>0\), and hence
\[
\frac{k_n}{n}\log r_n \to 0.
\]
Thus the bracketed term is \(-1+o(1)\). Since \(n/k_n\to\infty\),
\[
\log\!\left(
n\,\frac{\rho_n^{-n}r_n^\gamma}{k_n^t}
\right)
\le
-\frac{c}{2}\frac{n}{k_n}
\to -\infty,
\]
for all sufficiently large \(n\). Hence,
\[
\rho_n^{-n}
=
o\!\left(
\frac{k_n^t}{n}\,r_n^{-\gamma}
\right),
\]
which proves the result.
\end{proof}

\begin{lemma}[Lemma~\ref{lemma:lemma2_above_1_rho}]\label{lem:lemma4}
Under the mildly explosive regime described in Section~\ref{sec:section4},
\[
\frac{\rho_n^{-n}}{\ell_n k_n}
\sum_{t=1}^{n}\sum_{j=t}^{n}\rho_n^{\,t-j-1}u_j u_t
\;\xrightarrow{L^1}\;0,
\]
where
\[
\ell_n=\exp\!\left\{\frac{\alpha^2}{2(1-\phi_n^2)}\right\}.
\]
\end{lemma}

\begin{proof}
Decompose
\[
\frac{\rho_n^{-n}}{k_n\ell_n}
\sum_{t=1}^{n}\sum_{j=t}^{n}\rho_n^{\,t-j-1}u_j u_t
=
\frac{\rho_n^{-n-1}}{k_n\ell_n}\sum_{t=1}^{n}u_t^2
+
\frac{\rho_n^{-n}}{k_n\ell_n}
\sum_{t=1}^{n}\sum_{j=t+1}^{n}\rho_n^{\,t-j-1}u_j u_t
=:I_{1n}+I_{2n}.
\]
For the diagonal term,
\[
\mathbb E|I_{1n}|
=
\frac{\rho_n^{-n-1}}{k_n\ell_n}
\sum_{t=1}^n \mathbb E[u_t^2]
=
\frac{\rho_n^{-n-1}}{k_n\ell_n}
\sum_{t=1}^n x_t.
\]
Since \(x_t=\mathbb E[\sigma_t^2]\le C\ell_n\), we obtain
\[
\mathbb E|I_{1n}|
\le
C\frac{n\rho_n^{-n}}{k_n}
=o(1),
\]
by Lemma~\ref{lem:lemma3}. Hence \(I_{1n}\to0\) in \(L^1\).

\noindent For the off-diagonal term, expand the square. All off-diagonal terms vanish because the innovations \(\varepsilon_t\) are independent, centered, and independent of the volatility process. Thus only repeated pairs survive, and
\[
\mathbb E[I_{2n}^2]
=
\frac{\rho_n^{-2n}}{k_n^2\ell_n^2}
\sum_{t=1}^{n}
\sum_{j=t+1}^{n}
\rho_n^{2(t-j-1)}
\mathbb E[u_j^2u_t^2].
\]
By the lognormal moment bounds, there exists a finite constant \(q>0\) such that, uniformly in \(j,t\le n\),
\[
\mathbb E[u_j^2u_t^2]
=
\mathbb E[\sigma_j^2\sigma_t^2]
\le
C\ell_n^2 r_n^q.
\]
Therefore,
\[
\mathbb E[I_{2n}^2]
\le
C\frac{\rho_n^{-2n}r_n^q}{k_n^2}
\sum_{t=1}^{n}
\sum_{j=t+1}^{n}
\rho_n^{2(t-j-1)}.
\]
For each \(t\),
\[
\sum_{j=t+1}^{n}\rho_n^{2(t-j-1)}
\le
\sum_{h=1}^{\infty}\rho_n^{-2(h+1)}
=
\frac{\rho_n^{-2}}{\rho_n^2-1}
=
O(k_n),
\]
because \(\rho_n^2-1\sim 2c/k_n\). Hence
\[
\mathbb E[I_{2n}^2]
\le
C\frac{\rho_n^{-2n}r_n^q}{k_n^2}\, n k_n
=
C\frac{n\rho_n^{-2n}r_n^q}{k_n}.
\]
Now
\[
\frac{n\rho_n^{-2n}r_n^q}{k_n}
=
\rho_n^{-n}
\left(
\frac{n\rho_n^{-n}r_n^q}{k_n}
\right)
=o(1),
\]
again by Lemma~\ref{lem:lemma3}. Therefore
\[
\mathbb E[I_{2n}^2]\to0,
\]
so \(I_{2n}\to0\) in \(L^2\), and hence in \(L^1\). Combining \(I_{1n}\to0\) in \(L^1\) and \(I_{2n}\to0\) in \(L^1\) proves the claim.
\end{proof}

\begin{lemma}[Lemma~\ref{lemma:wn_vn_limit_rewrite_now}]\label{lemma:wn_vn_limit_rewrite}
Let $(W_n,V_n)$ be defined as in~\eqref{eq:def_wv_rewrite}. Then
\[
(W_n, V_n) \;\xRightarrow{d}\; (W, V) \quad \text{as } n \to \infty,
\]
where $W$ and $V$ are independent $\mathcal{N}\!\left(0,\frac{1}{2c}\right)$ random variables.
\end{lemma}
\begin{proof}
     We first verify the conditional variance condition. Recall
\[
\tau_{nj}
=
\frac{1}{\sqrt{\ell_n k_n}}
\left(
a\rho_n^{-j}
+
b\rho_n^{-(n-j+1)}
\right)u_j,
\qquad
u_j=\sigma_j\varepsilon_j .
\]
Then,
\[
aW_n+bV_n=\sum_{j=1}^n \tau_{nj}.
\]
We need to show that
\begin{align}
\sum_{j=1}^n
\mathbb E[\tau_{nj}^2\mid \mathcal F_{j-1}]
\xrightarrow{p}
\frac{a^2+b^2}{2c}.
\label{eq:cond-var-WV-target}
\end{align}
Since
\[
\mathbb E[u_j^2\mid\mathcal F_{j-1}]
=
\mathbb E[\sigma_j^2\varepsilon_j^2\mid \mathcal F_{j-1}]
=
\mathbb E[\sigma_j^2\mid z_{j-1}]
=
e^{\alpha^2/2}\sigma_{j-1}^{2\phi_n},
\]
we have
\begin{align}
\sum_{j=1}^n
\mathbb E[\tau_{nj}^2\mid\mathcal F_{j-1}]
=
\frac{e^{\alpha^2/2}}{\ell_n k_n}
\sum_{j=1}^n
\left(
a\rho_n^{-j}
+
b\rho_n^{-(n-j+1)}
\right)^2
\sigma_{j-1}^{2\phi_n}.
\label{eq:cond-var-WV}
\end{align}
Define
\[
\xi_{nj}
:=
\frac{e^{\alpha^2/2}}{\ell_n}\sigma_{j-1}^{2\phi_n}.
\]
Then \eqref{eq:cond-var-WV} becomes
\begin{align}
\sum_{j=1}^n
\mathbb E[\tau_{nj}^2\mid\mathcal F_{j-1}]
=
\frac{1}{k_n}
\sum_{j=1}^n
\left(
a\rho_n^{-j}
+
b\rho_n^{-(n-j+1)}
\right)^2
\xi_{nj}.
\label{eq:cond-var-xi}
\end{align}
It is easy to check, 
\begin{align}
\mathbb E[\xi_{nj}]
&=
\frac{e^{\alpha^2/2}}{\ell_n}
\exp\left\{
\frac{\alpha^2\phi_n^2(1-\phi_n^{2(j-1)})}
{2(1-\phi_n^2)}
\right\}
\nonumber\\
&=
\exp\left\{
-\frac{\alpha^2\phi_n^{2j}}
{2(1-\phi_n^2)}
\right\}.
\label{eq:mean-xi}
\end{align}
Hence, for \(j\ge M_n\),
\[
\mathbb E[\xi_{nj}]=1+o(1)
\]
uniformly, because
\[
\frac{\phi_n^{2M_n}}{1-\phi_n^2}=o(1).
\]
We now prove
\begin{align}
\frac{1}{k_n}\sum_{j=1}^n \rho_n^{-2j}\xi_{nj}
\xrightarrow{p}
\frac{1}{2c}.
\label{eq:left-weighted-lln}
\end{align}
By \eqref{eq:mean-xi},
\[
\begin{aligned}
\mathbb E\left[
\frac{1}{k_n}\sum_{j=1}^n \rho_n^{-2j}\xi_{nj}
\right]
&=
\frac{1}{k_n}\sum_{j=1}^{M_n}\rho_n^{-2j}\mathbb E[\xi_{nj}]
+
\frac{1}{k_n}\sum_{j=M_n+1}^{n}\rho_n^{-2j}\{1+o(1)\}.
\end{aligned}
\]
The first term is \(o(1)\), because \(M_n=o(k_n)\) and the moments are bounded
by a polynomial in \(r_n\), while \(k_n\) dominates such powers by Assumption~A3.
For the second term,
\[
\frac{1}{k_n}\sum_{j=M_n+1}^{n}\rho_n^{-2j}
=
\frac{\rho_n^{-2(M_n+1)}-\rho_n^{-2(n+1)}}
{k_n(\rho_n^2-1)}
\to
\frac{1}{2c},
\]
since \(M_n=o(k_n)\), \(\rho_n^{-2n}\to0\), and
\[
k_n(\rho_n^2-1)\to 2c.
\]
Thus,
\[
\mathbb E\left[
\frac{1}{k_n}\sum_{j=1}^n \rho_n^{-2j}\xi_{nj}
\right]
\to
\frac{1}{2c}.
\]
It remains to show concentration. For \(s<t\), set \(h=t-s\). From the Gaussian
log-volatility structure,
\[
\operatorname{Cov}(z_{s-1},z_{t-1})
=
\alpha^2\phi_n^h
\frac{1-\phi_n^{2(s-1)}}{1-\phi_n^2}
=
2\alpha^2\phi_n^h A_{s-1}.
\]
Hence
\[
\operatorname{Cov}(\xi_{ns},\xi_{nt})
=
\mathbb E[\xi_{ns}]\mathbb E[\xi_{nt}]
\left[
\exp\{2\alpha^2\phi_n^{h+2}A_{s-1}\}-1
\right].
\]
For \(h>M_n\),
\[
\phi_n^h A_{s-1}=O(\delta_n)=o(1),
\]
so
\[
|\operatorname{Cov}(\xi_{ns},\xi_{nt})|
\le C\delta_n.
\]
For \(h\le M_n\), the lognormal moment bounds give
\[
|\operatorname{Cov}(\xi_{ns},\xi_{nt})|
\le C r_n^q
\]
for some finite \(q>0\). Therefore
\[
\begin{aligned}
\operatorname{Var}
\left(
\frac{1}{k_n}\sum_{j=1}^n \rho_n^{-2j}\xi_{nj}
\right)
&\le
\frac{C}{k_n^2}
\sum_{j=1}^n \rho_n^{-4j} r_n^q
+
\frac{C}{k_n^2}
\sum_{h=1}^{M_n}
\sum_{j=1}^{n-h}
\rho_n^{-2j}\rho_n^{-2(j+h)}r_n^q
\\
&\quad+
\frac{C}{k_n^2}
\sum_{h>M_n}
\sum_{j=1}^{n-h}
\rho_n^{-2j}\rho_n^{-2(j+h)}\delta_n .
\end{aligned}
\]
Since,
\[
\sum_{j=1}^n \rho_n^{-2j}=O(k_n),
\]
the first two terms are bounded by
\[
C\frac{r_n^q}{k_n}
+
C\frac{M_n r_n^q}{k_n}
=o(1),
\]
by Assumption~A3. The last term is bounded by
\[
C\delta_n
\left(
\frac{1}{k_n}\sum_{j=1}^n \rho_n^{-2j}
\right)^2
=
O(\delta_n)
=o(1).
\]
Hence the variance converges to zero, proving \eqref{eq:left-weighted-lln}. Similarly,
\begin{align}
\frac{1}{k_n}
\sum_{j=1}^n
\rho_n^{-2(n-j+1)}
\xi_{nj}
\xrightarrow{p}
\frac{1}{2c}.
\label{eq:right-weighted-lln}
\end{align}
The proof is identical. The only change is the deterministic weight. Indeed,
\[
\frac{1}{k_n}
\sum_{j=M_n+1}^{n}
\rho_n^{-2(n-j+1)}
=
\frac{1-\rho_n^{-2(n-M_n)}}{k_n(\rho_n^2-1)}
\to
\frac{1}{2c},
\]
while the contribution from \(j\le M_n\) is negligible. Finally, the cross term is negligible. Since
\[
\rho_n^{-j}\rho_n^{-(n-j+1)}
=
\rho_n^{-(n+1)},
\]
we have
\[
\frac{1}{k_n}
\sum_{j=1}^n
\rho_n^{-j}\rho_n^{-(n-j+1)}
\xi_{nj}
=
\frac{\rho_n^{-(n+1)}}{k_n}
\sum_{j=1}^n \xi_{nj}.
\]
By Markov's inequality and \(\mathbb E[\xi_{nj}]\le C\),
\[
\mathbb E\left[
\left|
\frac{\rho_n^{-(n+1)}}{k_n}
\sum_{j=1}^n \xi_{nj}
\right|
\right]
\le
C\frac{n\rho_n^{-n}}{k_n}
\to0
\]
by Lemma~\ref{lem:lemma3}. Therefore
\begin{align}
\frac{\rho_n^{-(n+1)}}{k_n}
\sum_{j=1}^n \xi_{nj}
=o_p(1).
\label{eq:cross-weighted-negligible}
\end{align}
Combining \eqref{eq:left-weighted-lln}, \eqref{eq:right-weighted-lln}, and
\eqref{eq:cross-weighted-negligible} in \eqref{eq:cond-var-xi}, we obtain
\[
\sum_{j=1}^n
\mathbb E[\tau_{nj}^2\mid\mathcal F_{j-1}]
=
a^2\frac{1}{k_n}\sum_{j=1}^n\rho_n^{-2j}\xi_{nj}
+
b^2\frac{1}{k_n}\sum_{j=1}^n\rho_n^{-2(n-j+1)}\xi_{nj}
+
o_p(1),
\]
and hence
\[
\sum_{j=1}^n
\mathbb E[\tau_{nj}^2\mid\mathcal F_{j-1}]
\xrightarrow{p}
\frac{a^2+b^2}{2c}.
\]
This proves \eqref{eq:cond-var-WV-target}.

\noindent Now, we can verify the Lindeberg condition to prove asymptotic normality.

\noindent
We now verify the conditional Lindeberg condition. Recall that
\[
W_n=\frac{1}{\sqrt{\ell_n k_n}}\sum_{j=1}^n \rho_n^{-j}u_j,
\qquad
V_n=\frac{1}{\sqrt{\ell_n k_n}}\sum_{j=1}^n \rho_n^{-(n-j+1)}u_j.
\]
For fixed \(a,b\in\mathbb R\), define
\[
q_{nj}:=a\rho_n^{-j}+b\rho_n^{-(n-j+1)}
\]
and
\[
\tau_{nj}
=
\frac{q_{nj}}{\sqrt{\ell_n k_n}}u_j.
\]
Then
\[
aW_n+bV_n=\sum_{j=1}^n \tau_{nj}.
\]
For every fixed \(\eta>0\), define
\[
L_n(\eta)
:=
\sum_{j=1}^n
\mathbb E\left[
\tau_{nj}^2
\mathbf 1\{|\tau_{nj}|>\eta\}
\mid \mathcal F_{j-1}
\right].
\]
We show that
\begin{align}
L_n(\eta)=o_p(1).
\label{eq:lindeberg-ME-target}
\end{align}
we have
\begin{align}
L_n(\eta)
&\le
\frac{1}{\eta^2}
\sum_{j=1}^n
\mathbb E[\tau_{nj}^4\mid \mathcal F_{j-1}]
\nonumber\\
&=
\frac{1}{\eta^2\ell_n^2 k_n^2}
\sum_{j=1}^n
q_{nj}^4
\mathbb E[u_j^4\mid \mathcal F_{j-1}].
\label{eq:lindeberg-lyapunov-ME}
\end{align}
Since \(u_j=\sigma_j\varepsilon_j\), \(\mathbb E[\varepsilon_j^4]=\kappa_\varepsilon<\infty\), and
\[
z_j=\phi_n z_{j-1}+\eta_j,
\qquad
\eta_j\sim N(0,\alpha^2),
\]
we obtain
\[
\mathbb E[u_j^4\mid \mathcal F_{j-1}]
=
\kappa_\varepsilon
\mathbb E[\sigma_j^4\mid \mathcal F_{j-1}]
=
\kappa_\varepsilon
\mathbb E[e^{2z_j}\mid z_{j-1}].
\]
Therefore
\begin{align}
\mathbb E[u_j^4\mid \mathcal F_{j-1}]
=
\kappa_\varepsilon e^{2\alpha^2}\sigma_{j-1}^{4\phi_n}.
\label{eq:cond-fourth-ME}
\end{align}
Combining \eqref{eq:lindeberg-lyapunov-ME} and \eqref{eq:cond-fourth-ME}, it is enough to show
\begin{align}
\frac{1}{\ell_n^2 k_n^2}
\sum_{j=1}^n
q_{nj}^4\sigma_{j-1}^{4\phi_n}
=o_p(1).
\label{eq:ME-fourth-target}
\end{align}
We prove \eqref{eq:ME-fourth-target} by Markov's inequality. First, since
\[
q_{nj}
=
a\rho_n^{-j}+b\rho_n^{-(n-j+1)},
\]
there exists a constant \(C>0\), depending only on \(a\) and \(b\), such that
\[
q_{nj}^4
\le
C\left(
\rho_n^{-4j}
+
\rho_n^{-4(n-j+1)}
\right).
\]
Hence
\begin{align}
\mathbb E\left[
\frac{1}{\ell_n^2 k_n^2}
\sum_{j=1}^n
q_{nj}^4\sigma_{j-1}^{4\phi_n}
\right]
&\le
\frac{C}{\ell_n^2 k_n^2}
\sum_{j=1}^n
\left(
\rho_n^{-4j}
+
\rho_n^{-4(n-j+1)}
\right)
\mathbb E[\sigma_{j-1}^{4\phi_n}].
\label{eq:ME-fourth-expectation}
\end{align}
Now recall,
\[
\sigma_{j-1}^{4\phi_n}
=
e^{2\phi_n z_{j-1}}.
\]
Since
\[
z_{j-1}
\sim
N\left(
0,
\frac{\alpha^2(1-\phi_n^{2(j-1)})}{1-\phi_n^2}
\right),
\]
we get
\[
\mathbb E[\sigma_{j-1}^{4\phi_n}]
=
\exp\left\{
\frac{2\alpha^2\phi_n^2(1-\phi_n^{2(j-1)})}
{1-\phi_n^2}
\right\}
\le
\exp\left\{
\frac{2\alpha^2}{1-\phi_n^2}
\right\}.
\]
Since
\[
\ell_n
=
\exp\left\{
\frac{\alpha^2}{2(1-\phi_n^2)}
\right\},
\]
we have
\[
\ell_n^2
=
\exp\left\{
\frac{\alpha^2}{1-\phi_n^2}
\right\}.
\]
Thus
\begin{align}
\frac{\mathbb E[\sigma_{j-1}^{4\phi_n}]}{\ell_n^2}
&\le
\exp\left\{
\frac{\alpha^2}{1-\phi_n^2}
\right\}
\lesssim
r_n^{\alpha^2/(2d)}.
\label{eq:sigma-fourth-over-ell}
\end{align}
Using \eqref{eq:sigma-fourth-over-ell} in \eqref{eq:ME-fourth-expectation},
\[
\mathbb E\left[
\frac{1}{\ell_n^2 k_n^2}
\sum_{j=1}^n
q_{nj}^4\sigma_{j-1}^{4\phi_n}
\right]
\le
\frac{C r_n^{\alpha^2/(2d)}}{k_n^2}
\sum_{j=1}^n
\left(
\rho_n^{-4j}
+
\rho_n^{-4(n-j+1)}
\right).
\]
But
\[
\sum_{j=1}^n \rho_n^{-4j}
=
O(k_n),
\qquad
\sum_{j=1}^n \rho_n^{-4(n-j+1)}
=
O(k_n),
\]
because
\[
1-\rho_n^{-4}
\sim
\frac{4c}{k_n}.
\]
Therefore
\begin{align}
\mathbb E\left[
\frac{1}{\ell_n^2 k_n^2}
\sum_{j=1}^n
q_{nj}^4\sigma_{j-1}^{4\phi_n}
\right]
&\le
C\frac{r_n^{\alpha^2/(2d)}}{k_n}
=o(1),
\label{eq:ME-fourth-small}
\end{align}
by Assumption~\hyperref[assm:A3]{A3}. Hence, by Markov's inequality,
\[
\frac{1}{\ell_n^2 k_n^2}
\sum_{j=1}^n
q_{nj}^4\sigma_{j-1}^{4\phi_n}
=o_p(1).
\]
This proves \eqref{eq:ME-fourth-target}. Therefore, from
\eqref{eq:lindeberg-lyapunov-ME} and \eqref{eq:cond-fourth-ME},
\[
L_n(\eta)
=
\sum_{j=1}^n
\mathbb E\left[
\tau_{nj}^2
\mathbf 1\{|\tau_{nj}|>\eta\}
\mid \mathcal F_{j-1}
\right]
=o_p(1).
\]
Together with the conditional variance convergence,
\[
\sum_{j=1}^n
\mathbb E[\tau_{nj}^2\mid\mathcal F_{j-1}]
\xrightarrow{p}
\frac{a^2+b^2}{2c},
\]
the martingale central limit theorem implies
\[
aW_n+bV_n
=
\sum_{j=1}^n\tau_{nj}
\Rightarrow
N\left(0,\frac{a^2+b^2}{2c}\right).
\]
By the Cramér--Wold device,
\[
(W_n,V_n)\Rightarrow (W,V),
\]
where \(W\) and \(V\) are independent \(N(0,1/(2c))\) random variables.
Consequently, by the continuous mapping theorem and the decomposition established above,
\[
\left(
\frac{\rho_n^{-n}\sum_{t=1}^n y_{t-1}u_t}{\ell_n k_n},
\,
\frac{\rho_n^{-2n}\sum_{t=1}^n y_t^2}{\ell_n k_n^2}
\right)
\Rightarrow
(WV,W^2/2c).
\]
Thus the conditional Lindeberg condition holds.
\end{proof}

\noindent We now have all the tools needed to prove Theorem~\ref{thm:thm2}. Recall the time series process,
\begin{align*}
\frac{\rho_n^{-2n}}{k_n^2\,\ell_n}\sum_{t=1}^{n} y_{t-1}^2
&= \frac{1}{k_n^2\,\ell_n\,(\rho_n^2-1)}
\Bigg[
\rho_n^{-2n} (y_n^2 - y_0^2)
- 2\rho_n^{-2n+1}\sum_{t=1}^{n} y_{t-1}u_t
- \rho_n^{-2n}\sum_{t=1}^{n} u_t^2
\Bigg] \\[6pt]
&= \frac{1}{k_n\,\ell_n\,(\rho_n^2-1)}
\Bigg[
\frac{\rho_n^{-2n} y_n^2}{k_n}
- \frac{2\rho_n^{-2n+1}}{k_n}\sum_{t=1}^{n} y_{t-1}u_t
- \frac{\rho_n^{-2n}}{k_n}\sum_{t=1}^{n} u_t^2
\Bigg] + o(1)
\end{align*}
by Lemma~\ref{lemma:lemma1_app}. In addition, note that, by Lemma~\ref{lem:lemma3}
\[
\text{Var} \left(\frac{\rho_n^{-2n}}{k_n\,\ell_n}\sum_{t=1}^{n} u_t^2\right) = \mathcal{O}\!\,\left( \frac{n^2 \rho_n^{-4n} r_n^{\alpha^2/2d}}{k_n^2}\right) =o(1),
\]
\[
\frac{\rho_n^{-2n}}{k_n\,\ell_n}\sum_{t=1}^{n} u_t^2
= \mathcal{O}\!\,_p\!\left(\frac{n \rho_n^{-2n}r_n^{\frac{\alpha^2}{4d}}}{k_n}\right)
\,\,= o_p(1),
\]
and,
\[
\frac{\rho_n^{-2n+1}}{\ell_n k_n}
\sum_{t=1}^{n} y_{t-1}u_t
= \frac{\rho_n^{-n}y_0}{\sqrt{\ell_nk_n}}
\sum_{t=1}^{n} \frac{\rho_n^{-(n-t-1)} u_t}{\sqrt{\ell_n k_n}} + \frac{1}{\ell_n}\cdot\frac{\rho_n^{-2n+1}}{k_n}
\sum_{t=1}^{n}\left(\sum_{j=1}^{t-1}
\rho_n^{\,t-1-j}\,u_j\right)\,u_t.
\]
Also, note that the first term is $o(1)$ because of asymptotic normality of $V_n$, and the second term
\begin{align*}
\mathbb{E} \Bigg[
\frac{\rho_n^{-2n+1}}{k_n \ell_n}
\sum_{t=1}^{n}
\Bigg(
\sum_{j=1}^{t-1} \rho_n^{\,t-1-j} u_j
\Bigg) u_t
\Bigg]^2
&=
\frac{\rho_n^{-4n}}{\ell_n^2 k_n^2}
\sum_{t=1}^{n} \sum_{j=1}^{t-1}
\rho_n^{\,2(t-1-j)} \, E\Bigg[u_j^2u_t^2\Bigg] \\[6pt]
&\le
r_n^{2\alpha^2/d} \cdot
\frac{\rho_n^{-4n}}{k_n^2}
\sum_{t=1}^{n} \sum_{j=1}^{t-1}
\rho_n^{\,2(t-1-j)} \\[6pt]
&\le
\sum_{t=1}^{n}
(\rho_n^{\,2(t-1)}-1) r_n^{\alpha^2/2d}
\cdot
\frac{\rho_n^{-4n}}{\rho_n^2 - 1} \\[6pt]
&\le
\frac{n
\rho_n^{\,2n} r_n^{\alpha^2/2d} \rho_n^{-4n}
}{
k_n^2 (\rho_n^2 - 1)
} \\[6pt]
&= 
\mathcal{O}\!\,\!\left(\frac{n r_n^{\alpha^2/2d}
\rho_n^{-2n}}{k_n}
\right)
= o_p(1).
\end{align*}
Hence,
\begin{align}\label{eq:eq31}
\frac{\rho_n^{-2n}}{\ell_n k_n^2}
\sum_{t=1}^{n} y_{t-1}^2
&=
\frac{1}{k_n(\rho_n^2 - 1)}
\left(
\frac{\rho_n^{-n} y_n}{\sqrt{k_n \ell_n}}
\right)^2
+ o_p(1) \nonumber
\\[8pt]
&= \frac{1}{k_n(\rho_n^2 - 1)} \left(\frac{y_0}{\sqrt{k_n\ell_n}} + \sum_{j=1}^n \frac{\rho_n^{-j}u_j}{\sqrt{k_n\ell_n}}\right)^2 + o_p(1) \nonumber\\
&\xrightarrow{d}
\frac{1}{2c}\, W^2 ,
\qquad \text{where } W \sim N\!\left(0,\frac{1}{2c}\right).
\end{align}
This converges follows from the asymptotic normality of $W_n$.
Similarly, 
\begin{align*}
\frac{\rho_n^{-n}}{k_n \ell_n}
\sum_{t=1}^{n} y_{t-1} u_t
&=
\frac{\rho_n^{-n}}{k_n \ell_n}
\Bigg[
\sum_{t=1}^{n} y_0 \rho_n^{\,t-1} u_t
+
\sum_{t=1}^{n}
\Bigg(
\sum_{j=1}^{t-1} \rho_n^{\,t-1-j} u_j
\Bigg) u_t
\Bigg].
\end{align*}
Note,  from the asymptotic normality of $V_n$
\[
\sum_{t=1}^{n}
\frac{y_0}{\sqrt{k_n \ell_n}}
\rho_n^{-(\,n-t-1\,)}
\cdot
\frac{u_t}{\sqrt{k_n \ell_n}}
= o_p(1),
\]
Also using Lemma \ref{lem:lemma4} 
\begin{align}\label{eq:eq32}
\sum_{t=1}^{n}
\frac{\rho_n^{-n}}{k_n \ell_n}
\sum_{j=1}^{t-1}
\rho_n^{\,t-1-j} u_j u_t
=
\frac{\rho_n^{-n}}{k_n \ell_n}
\sum_{t=1}^{n}
\sum_{j=1}^{n}
\rho_n^{\,t-1-j} u_j u_t
- \frac{\rho_n^{-n}}{k_n \ell_n}
\sum_{t=1}^{n}
\sum_{j=t}^{n}
\rho_n^{\,t-1-j} u_j u_t, \nonumber \\
=W_nV_n +o_p(1) \xrightarrow{d}
WV \qquad \text{where} \, \,\, W, V\,\, \text{independent} \, \sim N(0, 1/2c)
\end{align}
Hence, combining~\eqref{eq:eq31} and~\eqref{eq:eq32} we conclude that
\begin{align}\label{eq:eq33}
\frac{
\dfrac{\rho_n^{-n} \sum_{t=1}^{n} y_{t-1}u_t}{\ell_n k_n}
}{
\dfrac{2c\rho_n^{-2n} \sum_{t=1}^{n} y_{t-1}^{2}}{\ell_n k_n^{2}}
}
&=
\frac{
\rho_n^{\,n} k_n \displaystyle\sum_{t=1}^{n} y_{t-1}u_t
}{
2c \displaystyle\sum_{t=1}^{n} y_{t-1}^{2}
} 
= \frac{
\rho_n^{\,n} k_n (\widehat{\rho_n} - \rho_n)}{
2c 
}
 = \frac{V}{W}\xrightarrow{d}
\mathcal{C}, 
\qquad
\mathcal{C} \sim \text{Cauchy}.
\end{align}
Thus, the proof of Theorem~\ref{thm:thm2} is complete.

\subsection{Proof of Theorem~\ref{thm:rho_gamma_ci_explosive}}
Consider the autoregression
\(
    y_t=\rho_n y_{t-1}+u_t,
    \,\,
    u_t=\sigma_t\varepsilon_t,
\)
under the regularity conditions of Theorem~\ref{thm:thm2}. 
Suppose $
    \rho_n=1+\frac{1}{n^{\gamma_n}},
    \,
    0<\gamma_n<1.$
Then
\begin{align}\label{main_rho_equation_CI_explosive}
    \frac{n^{\gamma_n}\rho_n^n}{2}
    \bigl(\widehat\rho_n-\rho_n\bigr)
    \Rightarrow \mathcal C,
\qquad 
\text{where}\,\, \mathcal C\sim\mathcal C(0,1).
\end{align} Moreover,
\begin{align}\label{main_rho_equation_CI_explosive}
    \rho_n^n\log n
    \bigl(\widehat\gamma_n-\gamma_n\bigr)
    \Rightarrow 2\mathcal C.
\end{align}
Consequently, if \(\widehat\rho_n>1\), an asymptotic \(100\lambda\%\) confidence interval for \(\rho_n\) is
\[
    \widehat\rho_n
    \pm
    z_\lambda^{\mathcal C}
    \frac{2}
    {n^{\widehat\gamma_n}\widehat\rho_n^{\,n}},
\qquad
    z_\lambda^{\mathcal C}
    =
    \Phi_{\mathcal C}^{-1}
    \!\left(\frac{1+\lambda}{2}\right)
\]
is the two-sided standard Cauchy critical value. The corresponding asymptotic \(100\lambda\%\) confidence interval for \(\gamma_n\) is
\[
    \widehat\gamma_n
    \pm
    z_\lambda^{\mathcal C}
    \frac{2}
    {
    \left(1+n^{-\widehat\gamma_n}\right)^n
    \log n
    }.
\]
\begin{proof}
We recall,
\begin{equation}\label{eq:exp_A_defs}
    A_n:=\rho_n-1=n^{-\gamma_n},
    \qquad
    \widehat A_n=\widehat\rho_n-1 .
\end{equation}
By Theorem~\ref{thm:thm2},
\begin{equation}\label{eq:exp_rho_cauchy}
    \frac{n^{\gamma_n}\rho_n^n}{2}
    (\widehat\rho_n-\rho_n)
    \xRightarrow{d}
    \mathcal C,
\end{equation}
where \(\mathcal C\sim\mathcal C(0,1)\). Equivalently, for some sequence
\(C_n\xRightarrow{d}\mathcal C\),
\begin{equation}\label{eq:exp_Cn_rep}
    \widehat\rho_n-\rho_n
    =
    \frac{2C_n}{n^{\gamma_n}\rho_n^n}.
\end{equation}
Using~\eqref{eq:exp_A_defs} and~\eqref{eq:exp_Cn_rep},
\begin{align}
    \frac{\widehat A_n}{A_n}
    &=
    \frac{A_n+(\widehat\rho_n-\rho_n)}{A_n} \notag\\
    &=
    1+n^{\gamma_n}(\widehat\rho_n-\rho_n) \notag\\
    &=
    1+\frac{2C_n}{\rho_n^n}.
    \label{eq:exp_A_ratio}
\end{align}
Since
\begin{equation}\label{eq:exp_rhon_growth}
    \rho_n^n
    =
    \left(1+n^{-\gamma_n}\right)^n
    =
    \exp\{n^{1-\gamma_n}(1+o(1))\}
    \to\infty,
\end{equation}
we have,
\begin{equation}\label{eq:exp_small_term}
    \frac{2C_n}{\rho_n^n}=o_p(1).
\end{equation}
Hence \(\widehat A_n/A_n=1+o_p(1)\). Moreover, from~\eqref{eq:exp_A_ratio},
\begin{equation}\label{eq:exp_abs_A}
    |\widehat A_n|
    =
    n^{-\gamma_n}
    \left|
    1+\frac{2C_n}{\rho_n^n}
    \right|.
\end{equation}
Recalling that,
\[
    \widehat\gamma_n
    =
    -\frac{\log|\widehat\rho_n-1|}{\log n}
    =
    -\frac{\log|\widehat A_n|}{\log n},
\]
and taking logarithms in~\eqref{eq:exp_abs_A}, we obtain
\begin{equation}\label{eq:exp_gamma_exact}
    \widehat\gamma_n-\gamma_n
    =
    -
    \frac{
    \log\left|
    1+\frac{2C_n}{\rho_n^n}
    \right|
    }{\log n}.
\end{equation}
By~\eqref{eq:exp_small_term}, the logarithmic perturbation is \(o_p(1)\). Therefore,
using \(\log(1+x)=x+o_p(x)\),
\begin{equation}\label{eq:exp_log_expansion}
    \log\left|
    1+\frac{2C_n}{\rho_n^n}
    \right|
    =
    \frac{2C_n}{\rho_n^n}
    +
    o_p\!\left(\rho_n^{-n}\right).
\end{equation}
Substituting~\eqref{eq:exp_log_expansion} into~\eqref{eq:exp_gamma_exact} gives
\begin{equation}\label{eq:exp_gamma_expansion}
    \widehat\gamma_n-\gamma_n
    =
    -
    \frac{2C_n}{\rho_n^n\log n}
    +
    o_p\!\left(
    \frac{1}{\rho_n^n\log n}
    \right).
\end{equation}
Consequently,
\begin{equation}\label{eq:exp_gamma_limit}
    \rho_n^n\log n
    (\widehat\gamma_n-\gamma_n)
    =
    -2C_n+o_p(1)
    \xRightarrow{d}
    2\mathcal C,
\end{equation}
where the last step uses the symmetry of the standard Cauchy distribution. It remains to justify the feasible normalizations. From~\eqref{eq:exp_gamma_limit},
\begin{equation}\label{eq:exp_gamma_log_rate}
    (\widehat\gamma_n-\gamma_n)\log n
    =
    O_p(\rho_n^{-n})
    =
    o_p(1).
\end{equation}
Therefore,
\begin{equation}\label{eq:exp_ngamma_ratio}
    \frac{n^{\widehat\gamma_n}}{n^{\gamma_n}}
    =
    \exp\{(\widehat\gamma_n-\gamma_n)\log n\}
    =
    1+o_p(1).
\end{equation}
Next, by~\eqref{eq:exp_Cn_rep},
\begin{equation}\label{eq:exp_rhohat_rate}
    \widehat\rho_n-\rho_n
    =
    O_p\!\left(
    \frac{1}{n^{\gamma_n}\rho_n^n}
    \right).
\end{equation}
Thus
\begin{align}
    n\log\!\left(\frac{\widehat\rho_n}{\rho_n}\right)
    &=
    n\log\!\left(
    1+\frac{\widehat\rho_n-\rho_n}{\rho_n}
    \right) \notag\\
    &=
    O_p\!\left(
    \frac{n}{n^{\gamma_n}\rho_n^n}
    \right)
    =
    O_p\!\left(
    \frac{n^{1-\gamma_n}}{\rho_n^n}
    \right)
    =
    o_p(1),
    \label{eq:exp_rhohat_log}
\end{align}
where the last equality follows from~\eqref{eq:exp_rhon_growth}. Hence
\begin{equation}\label{eq:exp_rhohat_power_ratio}
    \frac{\widehat\rho_n^{\,n}}{\rho_n^n}
    =
    \exp\!\left\{
    n\log\!\left(\frac{\widehat\rho_n}{\rho_n}\right)
    \right\}
    =
    1+o_p(1).
\end{equation}
Combining~\eqref{eq:exp_ngamma_ratio} and~\eqref{eq:exp_rhohat_power_ratio}, we get
\begin{equation}\label{eq:exp_feasible_delta_norm}
    n^{\widehat\gamma_n}\widehat\rho_n^{\,n}
    =
    n^{\gamma_n}\rho_n^n\{1+o_p(1)\}.
\end{equation}
Therefore, by~\eqref{eq:exp_rho_cauchy} and~\eqref{eq:exp_feasible_delta_norm},
\begin{equation}\label{eq:exp_feasible_delta_limit}
    \frac{
    n^{\widehat\gamma_n}\widehat\rho_n^{\,n}
    }{2}
    (\widehat\rho_n-\rho_n)
    \xRightarrow{d}
    \mathcal C,
\end{equation}
which yields the stated feasible confidence interval for \(\rho_n\).
Finally, from~\eqref{eq:exp_gamma_log_rate},
\[
    n^{-\widehat\gamma_n}
    =
    n^{-\gamma_n}\{1+o_p(1)\}.
\]
Using \(\rho_n=1+n^{-\gamma_n}\), it follows that
\begin{equation}\label{eq:exp_feasible_gamma_norm}
    \left(1+n^{-\widehat\gamma_n}\right)^n
    =
    \rho_n^n\{1+o_p(1)\}.
\end{equation}
Combining~\eqref{eq:exp_gamma_limit} with~\eqref{eq:exp_feasible_gamma_norm}, we obtain
\begin{equation}\label{eq:exp_feasible_gamma_limit}
    \left(1+n^{-\widehat\gamma_n}\right)^n
    \log n
    (\widehat\gamma_n-\gamma_n)
    \xRightarrow{d}
    2\mathcal C.
\end{equation}
This gives the stated feasible confidence interval for \(\gamma_n\).
\end{proof}

\end{spacing}    
\end{document}